# EDGEWORTH EXPANSIONS FOR STUDENTIZED STATISTICS UNDER WEAK DEPENDENCE[1]

By S. N. Lahiri

*Texas A & M University*

In this paper, we derive valid Edgeworth expansions for studentized versions of a large class of statistics when the data are generated by a strongly mixing process. Under dependence, the asymptotic variance of such a statistic is given by an infinite series of lag-covariances, and therefore, studentizing factors (i.e., estimators of the asymptotic standard error) typically involve an increasing number, say, $\ell$ of lag-covariance estimators, which are themselves quadratic functions of the observations. The unboundedness of the dimension $\ell$ of these quadratic functions makes the derivation and the form of the expansions *nonstandard*. It is shown that in contrast to the case of the studentized means under independence, the derived Edgeworth expansion is a *superposition of three distinct series*, respectively, given by one in powers of $n^{-1/2}$, one in powers of $[n/\ell]^{-1/2}$ (resulting from the standard error of the studentizing factor) and one in powers of the bias of the studentizing factor, where $n$ denotes the sample size.

**1. Introduction.** Studentized statistics (or t-statistics, in short) play a fundamental role in statistical inference about an unknown parameter. For example, construction of a confidence interval for a parameter $\theta$ and testing of hypotheses about $\theta$ are typically based on studentized statistics. However, in spite of its critical importance and uses, a complete understanding of the higher-order properties of studentized statistics for dependent observations remained elusive to date. As explained below in some details, a primary reason for this is that the standard Edgeworth expansion (EE) theory does not easily apply to the case of a studentized statistic under dependence. This paper develops the necessary tools and establishes valid EEs of a general order for a large class of studentized statistics under weak dependence. Such

Received March 2008; revised May 2009.
[1]Supported in part by NSF Grants DMS-07-42690 and DMS-07-07139.

EE results are important for investigating higher-order properties of tests and confidence intervals, such as coverage accuracy of one and two-sided block bootstrap confidence intervals and their iterated versions, power levels of tests under contiguous alternatives, etc., and for constructing inference procedures with higher-order accuracy.

To highlight the major theoretical issues associated with the problem, let $\{Z_i\}_{i=-\infty}^{\infty}$ be a sequence of stationary random variables with $EZ_1 = 0$ and $EZ_1^2 = \sigma^2 \in (0, \infty)$. Let $\bar{Z}_n = n^{-1} \sum_{i=1}^n Z_i$, $n \geq 1$. When the $Z_i$'s are independent and identically distributed (i.i.d.), $n^{1/2}\bar{Z}_n$ is asymptotically normal with mean zero and variance $\sigma^2$. In this case, $n^{1/2}\bar{Z}_n$ may be studentized using the scaling random variable

$$\hat{\sigma}_n^{\text{i.i.d.}} = \left[ n^{-1} \sum_{i=1}^n Z_i^2 - \bar{Z}_n^2 \right]^{1/2},$$

and therefore, the studentized sample mean

$$T_n^{\text{i.i.d.}} = n^{1/2}\bar{Z}_n / \hat{\sigma}_n^{\text{i.i.d.}}$$

may be expressed as a smooth function of the sample mean, $n^{-1}\sum_{i=1}^n (Z_i, Z_i^2)'$, of the bivariate random vectors $(Z_i, Z_i^2)', i = 1, \ldots, n$. Here and in the following, let $A'$ denote the transpose of a matrix $A$. An Edgeworth expansion for $T_n^{\text{i.i.d.}}$ can be derived using the well-known results on EEs for sums of i.i.d. random vectors and the transformation technique of Bhattacharya and Ghosh (1978).

The problem of deriving EEs for studentized sample mean becomes considerably more difficult when the $Z_i$'s are dependent. In the dependent case, under suitable moment and mixing conditions on the $Z_i$'s [cf. Ibragimov and Linnik (1971), Athreya and Lahiri (2006)], $n^{1/2}\bar{Z}_n$ is asymptotically normal with zero mean and variance

(1.1) $$\sigma_\infty^2 = \gamma(0) + 2 \sum_{k=1}^{\infty} \gamma(k),$$

where $\gamma(k) = \text{Cov}(Z_1, Z_{k+1})$, $k \in \mathbb{Z}$ and $\mathbb{Z} = \{0, \pm 1, \pm 2, \ldots\}$ is the set of all integers. Thus, for studentizing $\bar{Z}_n$, one needs to estimate all the lag-covariances $\gamma(k)$'s, eventually, as the sample size $n$ goes to infinity. In practice, one typically estimates a (large) number $\ell$ (say) of the lag-covariance terms and let $\ell = \ell_n$ tend to infinity with $n$ suitably. For example, with $\hat{\gamma}_n(k) \equiv n^{-1} \sum_{i=1}^{n-k} Z_i Z_{i+k} - \bar{Z}_n^2$, $1 \leq k < n$, a studentized version of $n^{1/2}\bar{Z}_n$ may be defined as

(1.2) $$T_n^{\text{DEP}} = n^{1/2}\bar{Z}_n / \hat{\sigma}_n^{\text{DEP}},$$

where $\hat{\sigma}_n^{\text{DEP}} = \max\{\hat{\gamma}_n(0) + 2\sum_{k=1}^{\ell} \hat{\gamma}_n(k), n^{-1}\}^{1/2}$ (or a suitable variant of it). However, note that unlike the i.i.d. case, the studentized statistic $T_n^{\text{DEP}}$ is now a function of $(\ell + 2)$ means, viz., $\bar{Z}_n$, $n^{-1} \sum_{i=1}^{n-k} Z_i Z_{i+k}, 0 \leq k \leq \ell$, where $\ell \to \infty$ as $n \to \infty$. It is precisely because of the *unbounded* dimensionality of the means in the definition of $T_n^{\text{DEP}}$ that the standard techniques of deriving EEs do not apply to the studentized statistics in the dependent case. Using some specialized arguments, two-term EEs [with errors of order $o(n^{-1/2})$] have been independently derived by Götze and Künsch (1996) and Lahiri (1996a). In a recent work, Velasco and Robinson (2001) obtained a third-order EE for the studentized sample mean for a *Gaussian* time series, using an exact Fourier inversion formula that is available only under the Gaussian assumption, but not in general. In this paper, we derive EEs of a *general order* for studentized versions of a large class of estimators based on different classes of commonly used studentizing factors under dependence. In particular, as regularity conditions on the underlying process, we only require some moment and weak dependence conditions as in Götze and Hipp (1983), but not Gaussianity. The arguments developed here builds upon a recent work of Lahiri (2007) and the seminal paper of Götze and Hipp (1983) (hereafter referred to as [L] and [GH], respectively) on EEs for sums of dependent variables and also, critically relies on a representation result of Bradley (1983).

To describe the main results of the paper, let $\{X_i\}_{i \in \mathbb{Z}}$ be a sequence of $\mathbb{R}^d$-valued stationary random vectors with $EX_1 = \mu$, where $d \in \mathbb{N} \equiv \{1, 2, \ldots\}$, the set of all positive integers. Let $\bar{X}_n = n^{-1} \sum_{i=1}^n X_i$ denote the sample mean of $X_1, \ldots, X_n$. Suppose that $\hat{\theta}_n$ is an estimator of a parameter of interest $\theta$ based on $X_1, \ldots, X_n$ such that

$$\hat{\theta}_n = H(\bar{X}_n) \quad \text{and} \quad \theta = H(\mu) \tag{1.3}$$

for some (smooth) function $H : \mathbb{R}^d \to \mathbb{R}$. This is the *smooth function* model of Bhattacharya and Ghosh (1978) [cf. Hall (1992)], which serves as a convenient, yet a general theoretical framework for studying EEs for a large class of commonly used statistics. Examples of statistics that can be treated under (1.3) include the sample mean, the sample auto-correlations and the Yule–Walker estimators of autoregressive parameters in an auto-regression model. The main results of the paper are also applicable to the class of generalized M-estimators of Bustos (1982), which are given by measurable solutions of estimating equations of the form

$$\sum_{i=1}^{n-p} \psi(X_i, \ldots, X_{i+p-1}; \hat{\theta}_n) = 0 \tag{1.4}$$

for some $\psi : \mathbb{R}^{dp+1} \to \mathbb{R}$, $p \in \mathbb{N}$ (see Remark 3.1 in Section 3 for details).

Note that under (1.3), the asymptotic variance of $n^{1/2}(\hat{\theta}_n - \theta)$ is given by

$$\tau_\infty^2 \equiv h(\mu)' \left[ \sum_{i=-\infty}^\infty \Gamma(i) \right] h(\mu),$$

where $\Gamma(i) = \text{Cov}(X_0, X_i) = E(X_0 - \mu)(X_i - \mu)'$, $i \in \mathbb{Z}$ and where $h(\cdot)$ is the $d \times 1$ vector of first-order partial derivatives of $H$. Here, we consider a class of studentizing factors of the form

$$\hat{\tau}_n^2 = \max\left\{ h(\bar{X}_n)' \left[ \hat{\Gamma}_n(0) + \sum_{k=1}^{\ell} w_{kn}\{\hat{\Gamma}_n(k) + \hat{\Gamma}_n(k)'\} \right] h(\bar{X}_n), n^{-1} \right\}$$

(1.5)
$$\equiv \max\{\hat{\tau}_{1n}^2, n^{-1}\} \quad \text{say}$$

(and also some of its commonly used variants; see Section 4), where $w_{kn} \in \mathbb{R}, 1 \leq k \leq \ell$ are nonrandom weights, where for $n \geq 1$, $\ell \equiv \ell_n \in (1, n)$, are integers such that $\ell_n^{-1} + n^{-1}\ell_n = o(1)$ as $n \to \infty$, and where

$$\hat{\Gamma}_n(k) \equiv n^{-1} \sum_{i=1}^{n-k} (X_i - \bar{X}_n)(X_{i+k} - \bar{X}_n)'$$

is a version of the sample auto-covariance matrix at lag $k$, $1 \leq k \leq n-1$. The studentized estimator is now given by

(1.6) $$T_n \equiv n^{1/2}(\hat{\theta}_n - \theta)/\hat{\tau}_n.$$

(For certain choices of $w_{kn}$'s, the first term in the definition of $\hat{\tau}_n^2$ can be negative or zero, albeit with small probability. The factor $n^{-1}$ in the definition of $\hat{\tau}_n^2$ ensures positivity of the scaling factor and thereby makes the studentized statistic $T_n$ well defined.) The main results of the paper give EEs for $T_n$ (and its variants) of order $s$ for any given integer $s \geq 3$. Not only are the arguments used for proving the EEs in the dependent case are different from those in the independent case, the forms of the EEs are also strikingly different. Recall that in the independent case, EEs for studentized statistics are given by a series in powers of $n^{-1/2}$:

(1.7) $$P(T_n^{\text{i.i.d.}} \leq x) = \Phi(x) + n^{-1/2} p_1(x)\phi(x) + n^{-1} p_2(x)\phi(x) + \cdots,$$

uniformly in $x \in \mathbb{R}$, where $\Phi(\cdot)$ and $\phi(\cdot)$, respectively, denote the distribution and the density functions of a $N(0,1)$ variate, and $p_i(\cdot)$'s are some polynomials. In contrast, *the EEs in the dependent case are super-impositions of three distinct series in powers of $n^{-1/2}$, $b_n^{-1/2}$ and $a_n^{-1}$* [cf. (3.8), Proposition 3.2 below]:

(1.8)
$$P(T_n \leq x) = \Phi(x) + a_n^{-1} q_1(x)\phi(x) + n^{-1/2} q_2(x)\phi(x)$$
$$+ b_n^{-1} q_3(x)\phi(x) + \cdots,$$

where $q_i(\cdot)$'s are some polynomials, $b_n = n/\ell$ and $a_n^{-1}$ is of the same order as the bias of $\hat{\tau}_n^2$. In (1.8), the $n^{-1/2}$-series results from the centered and scaled sample mean $n^{1/2}(\bar{X}_n - \mu)$ as in [GH]. The terms in the $a_n^{-1}$-series and the $b_n^{-1/2}$-series result from the estimation of the asymptotic variance $\tau_\infty^2$ by the "truncated series" estimator $\hat{\tau}_n^2$; Here, the $b_n^{-1/2}$-series is determined by the stochastic part of $\hat{\tau}_n^2$ and the $a_n^{-1}$-series is determined by the deterministic truncation error (or the bias of the estimator). We point out that the coefficient of the $b_n^{-1/2}$-term in the EE of $T_n$ in (1.8) is zero, but those of the higher powers of $b_n^{-1/2}$ typically are not. To provide further insight into the interactions among the three series and for ready reference, we provide explicit expressions for the terms in the EE for $T_n$ with an error of order $o(n^{-1})$ in Section 3 below. It is evident from this explicit EE [cf. (3.8)] that interactions between the $n^{-1/2}$-series and the $b_n^{-1/2}$-series show up *only* in the terms of order *smaller* than $n^{-1/2}$, and hence, the existing second-order EE results fail to provide any clue to the complex pattern of interactions between these series. In addition to the general-order EE for $T_n$, we also establish valid EEs of a general order under alternative forms of studentizations, including those based on various block bootstrap and subsampling methods.

The rest of the paper is organized as follows. We conclude this section with a brief literature review. In Section 2, we introduce the theoretical frameworks of [GH] and [L] and state the assumptions. Here, we also verify these assumptions for a class of vector linear processes. In Section 3, we present the main results under the smooth function model. EE results under other commonly used choices of the studentizing factors are given in Section 4. Proofs of the main results are given in Sections 5 and 6.

There is a vast literature on EEs for the sample mean $\bar{X}_n$ and for statistics that are smooth functions of $\bar{X}_n$. For *independent* random vectors, a detailed account of the EE theory for $\bar{X}_n$ is given by Bhattacharya and Ranga Rao (1986) (hereafter referred to as [BR]) and Petrov (1975), and the theory under the "smooth function model" is given by Bhattacharya and Ghosh (1978), Bhattacharya (1985), and Hall (1992), among others. For *weakly dependent* random vectors, [GH] obtained EEs for $\bar{X}_n$ under a very flexible framework. Lahiri (1993) relaxed the moment condition used by [GH] and settled a conjecture of [GH] on the validity of expansions for expectations of smooth functions of $\bar{X}_n$. Applicability of [GH] results in different time series models have been verified in Götze and Hipp (1994). EEs for $\bar{X}_n$ under polynomial mixing rates have been given by Lahiri (1996b). Second-order expansions for certain versions of studentized statistics under weak dependence are given by Götze and Künsch (1996) and Lahiri (1996a). Velasco and Robinson (2001) derived third-order edgeworth expansion for a studentized version of the sample mean for a stationary Gaussian process. EEs for sums of block-variables are recently given in [L]. The results of this paper establish EEs for most commonly used versions of studentized statistics under weak dependence and extend earlier results on the problem,

6  S. N. LAHIRIwhere the complex pattern of interactions among the three series were not very clear.

**2. Theoretical framework and assumptions.** For deriving EEs for the studentized statistic $T_n$, we will adopt a general framework similar to the one introduced by [GH] in the context of establishing valid asymptotic expansions for sums of weakly dependent random vectors. Suppose that the stationary sequence of random vectors $\{X_i\}_{i\in\mathbb{Z}}$ are defined on a probability space $(\Omega, \mathcal{F}, P)$. Also, suppose that $\{\mathcal{D}_j\}_{j=-\infty}^{\infty}$ is a given collection of sub $\sigma$-fields of $\mathcal{F}$. Let $\mathcal{D}_p^q \equiv \sigma\langle \mathcal{D}_j : j \in \mathbb{Z}, p \leq j \leq q\rangle, -\infty \leq p \leq q \leq \infty$. For any two sub-sigma-fields $\mathcal{G}, \mathcal{H}$ of $\mathcal{F}$, define the $\alpha$-mixing coefficient between $\mathcal{G}$ and $\mathcal{H}$ as $\alpha(\mathcal{G}, \mathcal{H}) = \sup\{|P(A \cap B) - P(A)P(B)| : A \in \mathcal{G}, B \in \mathcal{H}\}$. *For notational simplicity, we set $\mu = EX_1 = 0$ in the statements of the assumptions below, unless there is some scope of confusion.* [If in an application $\mu \neq 0$, then replace $X_i$'s in (A.1)–(A.6) below with $(X_i - \mu)$'s.] For $i \in \mathbb{Z}$, let $X_i = (X_{i,1}, \ldots, X_{i,d})'$ and let $Y_{in}^{\#}$ be a $[d_1 \equiv d(d+1)/2]$-dimensional vector with $(p,q)$th element

$$Y_{in}^{\#}(p,q) = X_{i,p}X_{i,q} + \sum_{k=1}^{\ell} w_{ikn}[X_{i,p}X_{i+k,q} + X_{i+k,p}X_{i,q}],$$

(2.1)
$$1 \leq p \leq q \leq d,$$

where $w_{ikn} = w_{kn}$ for $1 \leq i + k \leq n$ and $w_{ikn} = 0$ for $i + k > n$. Next, define the *block variables*

$$(2.2) \quad W_{jn} = \left(\ell^{-1/2} \sum_{i=(j-1)\ell+1}^{j\ell \wedge n} X_i'; \ell^{-1} \sum_{i=(j-1)\ell+1}^{j\ell \wedge n} Y_{in}'\right)', \qquad 1 \leq j \leq b_{0n},$$

where $Y_{in} = Y_{in}^{\#} - EY_{in}^{\#}$ and $b_{0n} = \lceil n/\ell \rceil$. Here and in the following, for $x \in \mathbb{R}$, we write $\lceil x \rceil$ to denote the smallest integer not less than $x$ and $\lfloor x \rfloor$ to denote the integer part of $x$. Also, we write $\|\cdot\|$ to denote the Euclidean norm of a vector and the spectral norm of a matrix, i.e., $\|(x_1, \ldots, x_k)'\| = [\sum_{i=1}^{k} x_i^2]^{1/2}$, $(x_1, \ldots, x_k)' \in \mathbb{R}^k$ and $\|A\| = \sup\{\|Ax\| : x \in \mathbb{R}^k, \|x\| = 1\}$ for a $r \times k$ matrix $A$, where $r, k \in \mathbb{N} \equiv \{1, 2, \ldots\}$. For a random vector $U$ and $\gamma \in (0, \infty)$, let $\|U\|_{\gamma} = (E\|U\|^{\gamma})^{1/\gamma}$. Let $f \equiv ((f_{ij}))$ denote the $(d \times d$ matrix-valued) spectral density function of the process $\{X_i\}_{i\in\mathbb{Z}}$. Set $\iota = \sqrt{-1}$ and $\delta_n = b_n^{-(s-2)/2}(\log n)^{-2}$, $n \geq 2$.

We shall make use of the following assumptions for deriving an $(s-2)$th "order" EE for $T_n$, for a given integer $s \geq 3$. Note that an $(s-2)$th "order" EE here is interpreted as a uniform approximation to the distribution of $T_n$ with an *error term* $o(b_n^{-(s-2)/2})$. This may be contrasted with the more familiar notion of order in classical EEs in the independent case (e.g., for $T_n^{\text{i.i.d.}}$), where the error of an $(s-2)$th-order EE is $o(n^{-(s-2)/2})$. This distinction is important as the error of approximation in the dependent case is determined by a coarser series (in powers of $b_n^{-1/2}$ compared to a series in powers of $n^{-1/2}$), resulting from estimation of the infinite-dimensional parameter $f(0)$.

ASSUMPTIONS.

(A.1) (i) There exists a constant $\kappa \in (0,1)$ such that for all $n \geq \kappa^{-1}$, $\kappa \log n < \ell < \kappa^{-1} n^{1/2-\kappa}$.

(ii) The function $H : \mathbb{R}^d \to \mathbb{R}$ in (1.3) is $[\nu(s)+2]$-times differentiable in a neighborhood of $\mu = EX_1$ and $h(\mu) \neq 0$, where

$$\nu(s) = \inf\{r \in \mathbb{N} : n^{-r/2}(\log n)^{(r+1)/2} = O(\delta_n) \text{ as } n \to \infty\}$$

and where recall that $h(\cdot)$ is the $d \times 1$ vector of first-order partial derivatives of $H$.

(A.2) (i) There exists a constant $\kappa \in (0,1)$ such that $\sup\{\|W_{jn}\|_{s+\kappa} : 1 \leq j \leq b_{0n}, n \geq 1\} < \kappa^{-1}$.

(ii) The $d \times d$ matrix $f(0)$ is nonsingular and $\Delta_\infty \equiv \lim_{n\to\infty} \text{Cov}(\frac{1}{\sqrt{n\ell}} \times \sum_{i=1}^n Y_{in})$ exists and is nonsingular.

(iii) $\sup\{|w_{kn}| : 1 \leq k \leq \ell, n \geq 1\} < \infty$ and $\lim_{n\to\infty} |w_{kn} - 1| = 0$ for each $k \geq 1$.

(A.3) There exists a constant $\kappa \in (0,1)$ such that for all $m > \kappa^{-1}$ and for all $j \geq 1$, there exist a $\mathcal{D}_{j-m}^{j+m}$-measurable $X_{j,m}^\dagger$ such that

$$E\|X_j - X_{j,m}^\dagger\| \leq \kappa^{-1} \exp(-\kappa m).$$

(A.4) There exists a constant $\kappa \in (0,1)$ such that for all $n, m = 1, 2, \ldots$,

$$\alpha(\mathcal{D}_{-\infty}^n, \mathcal{D}_{n+m}^\infty) \leq \kappa^{-1} \exp(-\kappa m).$$

(A.5) There exists a constant $\kappa \in (0,1)$ such that for all $i, j, k, r, m = 1, 2, \ldots$ with $i < k < r < j$ and $m > \kappa^{-1}$,

$$E\Big(\sup_{A \in \mathcal{D}_i^j} |P(A|\mathcal{D}_j : j \notin [k,r]) - P(A|\mathcal{D}_j : j \in [i-m,k) \cup (r, j+m])|\Big)$$
$$\leq \kappa^{-1} \exp(-\kappa m).$$

(A.6) There exists a constant $\kappa \in (0,1)$ and a sequence $\{d_n\} \subset [1, \infty)$ with $d_n = O(\ell)$ and $d_n^2 = O(b_n^{1-\kappa})$ such that for all $n \geq \kappa^{-1}$ and $\kappa^{-1} < m \leq [\log n]^{1/\kappa}$,

(2.3) $$\max_{j_0 \in J_n} E\Bigg[\sup_{t \in B_n} \bigg| E\bigg\{\exp\bigg(\iota t' \sum_{j=j_0-m}^{j_0+m} W_{jn}\bigg) \bigg| \tilde{\mathcal{D}}_{j_0}\bigg\} \bigg|\Bigg] \leq 1 - \kappa,$$

where $J_n \equiv J_n(m) = \{m+1,\ldots,b_{0n}-m-1\}$, $B_n = \{t \in \mathbb{R}^d : \kappa d_n \leq \|t\| \leq b_n^{(s-2)/2+\kappa}\}$, and $\tilde{\mathcal{D}}_{j_0} \equiv \tilde{\mathcal{D}}_{j_0}(m,\ell) = \sigma\langle\{\mathcal{D}_j : j \in \mathbb{Z}, j \notin [(j_0 - \lfloor \frac{m}{2} \rfloor)\ell + 1, (j_0 + \lfloor \frac{m}{2} \rfloor + 1)\ell]\}\rangle$.

Now we comment on the assumptions. As in [GH], here we formulate the assumptions in terms of the auxiliary $\sigma$-fields $\mathcal{D}_j$'s in order to allow for more flexibility and generality in applications. (A.1)(i) is a growth condition on the block length $\ell$ and allows $\ell$ to grow at a rate $O(n^{[1/2]-\kappa})$ for an arbitrarily small $\kappa > 0$. Most of the commonly used variance estimators use an $\ell$ in the range specified by (A.1)(i). For example, for the overlapping or nonoverlapping block bootstrap variance estimators based on blocks of length $\ell$, the (MSE-) optimal rate for $\ell$ is $O(n^{1/3})$ [cf. Künsch (1989), Hall, Horowitz and Jing (1995), Lahiri (2003)], for weights $w_{kn}$'s based on symmetric kernels, the optimal rate of $\ell$ is $O(n^{1/5})$ [cf. Priestley (1981)], and for weights based on flat-top kernels, the optimal rate of $\ell$ is $O(\log n)$ under exponential strong mixing conditions [cf. Politis and Romano (1995)]. (A.1)(ii) is a smoothness assumption on the function $H$. The order $\nu(s)$ is determined by the relative sizes of $n$ and $\ell$. In particular, if $H$ is $(s+1)$-times differentiable in a neighborhood of $\mu$ with $h(\mu) \neq 0$, then (A.1)(ii) holds for all choices of $\ell$ satisfying (A.1)(i).

Assumption (A.2)(i) is a moment condition on the block variables $W_{jn}$'s and can be weakened slightly, where boundedness of the $(s+\kappa)$th absolute moments of $W_{jn}$'s is replaced by $\sup\{E\|W_{jn}\|^s \log(1+\|W_{jn}\|)^{\alpha(s)} : j = 1,\ldots,b_{0n}, n \geq 1\} < \kappa^{-1}$ for some $\alpha(s) > 2s^2$ and $\kappa \in (0,1)$ (cf. [L]). A simple sufficient condition for (A.2)(i), assuming (A.2)(iii), is that $E\|X_1\|^{2(s+1)+\kappa} < \infty$ for some $\kappa > 0$. This can be proved using the arguments used in the proof of Lemma 3.28 of [GH].

Next consider assumption (A.2)(ii). The two parts of (A.2)(ii) jointly imply that the limiting covariance matrix of the normalized sum $b_{0n}^{-1/2} \sum_{j=1}^{b_{0n}} W_{jn}$ exists and is *nonsingular*. Note that when the weights $w_{kn}$'s are generated through a bounded weight function $w:[0,1] \to \mathbb{R}$ as

$$(2.4) \qquad w_{kn} = w(k/\ell), \qquad 0 \leq k \leq \ell,$$

simple sufficient conditions on $w(\cdot)$ guaranteeing the existence of $\Delta_\infty$ are known [cf. page 699, Priestley (1981)]. Also, in this case, (A.2)(iii) holds if $w(0) = 1$. However, the main results of the paper remain valid for more general weights, as long as assumption (A.2) holds, and therefore, it is not necessary to assume that the weights be generated through (2.4).

Assumption (A.3) is an approximation condition that connects the variables $X_i$'s to the strong-mixing property (A.4) of the auxiliary $\sigma$-fields $\mathcal{D}_j$'s. The exponential decay of the strong mixing co-efficient in (A.4) can be weakened to a suitable polynomial rate as in Lahiri (1996b), but at the expense of a much longer proof. We do not pursue such refinements here to save space. Assumption (A.5) is an approximate Markov condition and is a variant of a similar condition used by [GH]. This condition holds if the original random vectors $X_i$'s are generated by a Markov chain. Although assumptions (A.3)–(A.5) may seem to have some overlap, in general, none of (A.3)–(A.5) implies the other(s). However, in a specific application, it may be possible to deduce one or more of these assumptions from the rest and can be dropped. For example, if the random vectors $\{X_i\}_{i=-\infty}^{\infty}$ are strongly mixing at an exponential rate and if we set $\mathcal{D}_i = \sigma\langle X_i \rangle$ for all $i$, then assumption (A.4) holds and (A.3) becomes redundant (by choosing $X_{i,m} = X_i$ for all $m \geq 1$).

Finally, consider assumption (A.6), which is a stronger version of a conditional Cramer's condition on the block variables $W_{jn}$'s given in [L] (see also condition (2.6) in [GH]). From the proof of the main result (Theorem 3.1), it is evident that the studentized statistic under consideration is a smooth function of a variable of the form $b_{0n}^{-1/2} \sum_{j=1}^{b_{0n}} W_{jn} + b_n^{-1} \xi_{1n}$, where the components of $\xi_{1n}$ are certain quadratic functions of central and boundary block variables. To handle the special structure of $\xi_{1n}$, we employ a representation theorem of Bradley (1983) for strongly mixing random variables and develop some iterated conditioning arguments to separate out the effects of the boundary block variables from the central ones. It is worth noting that typically, other versions of mean-corrected studentized statistic have a similar structure and, therefore, the arguments developed here may be of some independent interest in the context of studentization under weak dependence. The stronger version of the standard conditional Cramer's condition (i.e., the uniformity requirement over $t \in B_n$) is a consequence of the iterated conditioning step in the proof.

Verification of assumptions (A.3)–(A.6) requires choosing the auxiliary $\sigma$-fields $\mathcal{D}_j$'s and the sequence $d_n$ suitably. Typically, the naive choice $\mathcal{D}_i = \sigma\langle X_i \rangle$ is not the best to work with. As an illustration, we now consider the important special case where $\{X_i\}_{i \in \mathbb{Z}}$ is a vector linear process and provide some simple sufficient conditions for (A.3)–(A.6).

PROPOSITION 2.1. *Let*

$$(2.5) \qquad X_i = \mu + \sum_{j \in \mathbb{Z}} A_j \varepsilon_{i-j}, \qquad i \in \mathbb{Z},$$

*where $\{A_i : i \in \mathbb{Z}\}$ is a collection of nonrandom $d \times d$ matrices and $\{\varepsilon_i\}_{i \in \mathbb{Z}}$ is a sequence of $\mathbb{R}^d$-valued i.i.d. random vectors with $E\varepsilon_1 = 0$. Suppose that $\sup\{|w_{kn}| : 1 \leq k \leq \ell, n \geq 1\} < \infty$, $\lim_{n \to \infty} w_{kn} = 1$ for each $k \geq 1$, $\|A_i\| = O(\exp(-\kappa|i|))$ as $|i| \to \infty$, for some $\kappa \in (0,1)$ and that $A_\infty \equiv \sum_{i \in \mathbb{Z}} A_i$ is nonsingular. Further, suppose that $\varepsilon_1$ has a nonzero absolutely continuous component with respect to the Lebesgue measure on $\mathbb{R}^d$. Then, assumptions (A.3)–(A.6) hold.*

Now, we indicate how (A.3)–(A.6) are verified for the linear process case above. Under (2.5), we set $\mathcal{D}_j = \sigma\langle\varepsilon_j\rangle$, $j \in \mathbb{Z}$. By the independence of the $\varepsilon_j$'s, assumptions (A.4) and (A.5) hold trivially. For (A.3), we take $X^\dagger_{i,m} = \mu + \sum_{|j|\leq m} A_j \varepsilon_{i-j}$, $i \in \mathbb{Z}, m \geq 1$. Then, $X^\dagger_{i,m}$ is $\mathcal{D}^{i+m}_{i-m}$-measurable and in view of the exponential decay of $\|A_i\|$, (A.3) holds. Verification of (A.6) requires some additional work; see the proof of Proposition 2.1 in Section 6 for details.

**3. Main results.** For describing the form of the EE and for stating the main results, we will introduce some notation and conventions at this point. Let $\mathbb{Z}_+ = \{0, 1, 2, \ldots\}$. For $\alpha = (\alpha_1, \ldots, \alpha_d)' \in \mathbb{Z}_+^d$ and $x = (x_1, \ldots, x_d) \in \mathbb{R}^d$, set $|\alpha| = \alpha_1 + \cdots + \alpha_d$, $\alpha! = \prod_{k=1}^d \alpha_k$, and $x^\alpha = \prod_{k=1}^d x_k^{\alpha_k}$, and let $D^\alpha$ denote the differential operator $\frac{\partial^{\alpha_1}}{\partial x_1^{\alpha_1}} \cdots \frac{\partial^{\alpha_d}}{\partial x_d^{\alpha_d}}$ where $\frac{\partial^a}{\partial x_i^a}$ denotes $a$th-order partial derivative with respect to the $i$th coordinate of $x$. Write $|I|$ to denote the size of a finite set $I$. For a collection of variables $\{A(\lambda) : \lambda \in \Lambda\}$ and for $I \subset \Lambda$, set $A(I) = \prod_{\lambda \in I} A(\lambda)$ if $I \neq \varnothing$ and $A(I) = 1$ if $I = \varnothing$. Unless otherwise stated, (i) limits in the order symbols are taken by letting $n \to \infty$ and, (ii) the components of a vector indexed by the elements of the set $\Lambda_0 \equiv \{(p,q) : p, q \in \mathbb{N}, 1 \leq p \leq q \leq d\}$ follow the *lexicographic ordering*. Let $\mathbb{1}(\cdot)$ denote the indicator function.

3.1. *Definition of the EE.* Using the smoothness of the function $H$, we can use Taylor's expansion to derive a stochastic expansion $T_{1n}$ (say) for $T_n$ such that $T_{1n}$ is a polynomial function of certain vector-valued block-variables and such that the $(s-2)$th-order EEs of $T_n$ and $T_{1n}$ match up to an error of $o(b_n^{-(s-2)/2})$. In Section 6, we show that (a choice of) $T_{1n}$ is given by

$$T_{1n} = \frac{h(\mu)' Z_n}{\tau_{1n}} + \sum_{j=1}^{\nu(s)} n^{-j/2} \sum_{|\alpha|=j+1} c_{1n}(\alpha) Z_n^\alpha$$

(3.1)

$$+ \sum_{k=1}^{s-2} \sum_{j=0}^{\nu_k(s)} n^{-j/2} b_n^{-k/2} \sum_{|\alpha|=j+1} \sum^{*k} c_{2n}(\alpha, I) \prod_{r=1}^3 Z_n^\alpha \hat{A}_{rn}(I_r),$$

where $Z_n \equiv n^{1/2}(\bar{X}_n - \mu)$, $\tau_{1n}^2 = h(\mu)'[\Gamma(0) + \sum_{k=1}^\ell w_{kn}\{\Gamma(k) + \Gamma(k)'\}]h(\mu)$, $\nu_k(s) = \min\{r \in \mathbb{Z}_+ : n^{-(r+1)/2} b_n^{-k/2}(\log n)^{(k+r)/2} = O(\delta_n)\}$, and where $\{c_{1n}(\alpha)\}$ and $\{c_{2n}(\alpha, I)\}$ are some *bounded* sequences of real numbers, $\hat{A}_{kn} = (\hat{A}_{kn}(\lambda) : \lambda \in \Lambda_0)$, $k = 1, 2, 3$ are certain $d_1 \times 1$ variables, as defined in (6.46), and satisfy $\|\hat{A}_{kn}\| = O_p(1)$. The sum $\sum^{*k}$ in (3.1) extends over all $I = (I_1, I_2, I_3)$ such that $I_j \subset \Lambda_0$, $j = 1, 2, 3$ and $|I_1| + 2|I_2| + 3|I_3| = k$, $1 \leq k \leq s-2$. Let $\mathcal{X}_{r,n}$ denote the $r$th cumulant of $T_{1n}$, $r \in \mathbb{N}$ (when it exists) and let

$\tau_n^2 = \text{Var}(h(\mu)'Z_n)$. Next, we expand the cumulants of $T_{1n}$ *formally* (i.e., assuming existence of all moments of the variables appearing in $T_{1n}$) to get

(3.2)
$$\mathcal{X}_{r,n} = \sum_{k=0}^{s-2}\sum_{j=0}^{\nu_k(s)} b_n^{-k/2} n^{-j/2} \beta_{k,j}^{(r)} + o(b_n^{-(s-2)/2}) \quad \text{for } 1 \leq r \leq s, r \neq 2,$$

$$\mathcal{X}_{2,n} = \frac{\tau_n^2}{\tau_{1n}^2} + \sum_{k=0}^{s-2}\sum_{j=0}^{\nu_k(s)} b_n^{-k/2} n^{-j/2} \beta_{k,j}^{(2)} + o(b_n^{-(s-2)/2}) \quad \text{for } r = 2,$$

where $\beta_{k,j}^{(r)} \equiv \beta_{k,j,n}^{(r)}$ are real numbers satisfying $\beta_{k,j,n}^{(r)} = O(1)$ for each $j, k$, with $\beta_{0,0,n}^{(r)} \equiv 0$ for all $n \geq 1$ and $r \in \{1, \ldots, s\}$, and where $\nu_0(s) = \nu(s) - 2$. Then define the polynomials $\pi_{j,k,n}(\cdot)$ by the identity (in $t \in \mathbb{R}$)

$$\exp\left(\sum_{r=1}^{s}\left[\sum_{k=0}^{s-2}\sum_{j=0}^{\nu_k(s)} b_n^{-k/2} n^{-j/2} \beta_{k,j}^{(r)}\right](\iota t)^r / r!\right)$$

$$= 1 + \sum_{k=0}^{s-2}\sum_{j=0}^{\nu_k(s)} b_n^{-k/2} n^{-j/2} \pi_{j,k,n}(\iota t) + o(b_n^{-(s-2)/2}),$$

where $\pi_{0,0,n}(\cdot) \equiv 0$ and for $(j,k) \neq (0,0)$, the coefficients of $\pi_{j,k,n}$ are functions of $\beta_{k,j}^{(r)} = \beta_{k,j,n}^{(r)}$'s and are *bounded* [i.e., $O(1)$ as $n \to \infty$]. (We show the subscript $n$ in $\pi_{j,k,n}$ to highlight the dependence on $\beta_{k,j,n}^{(r)}$'s.) The density $\psi_{s,n}^*$ of the $(s-2)$*th-order* preliminary EE for $T_n$ (with respect to the Lebesgue measure on $\mathbb{R}$) is given by

(3.3) $$\psi_{s,n}^*(u) = \left[1 + \sum_{k=0}^{s-2}\sum_{j=0}^{\nu_k(s)} b_n^{-k/2} n^{-j/2} \pi_{j,k,n}\left(-\frac{d}{du}\right)\right]\phi_{e_n}(u), \quad u \in \mathbb{R},$$

where $\pi_{0,0,n}(\cdot) \equiv 0$, $e_n^2 = \tau_n^2/\tau_{1n}^2$ and $\phi_a(u) = a^{-1}\phi(u/a)$, $u \in \mathbb{R}, a \in (0, \infty)$. Note that in (3.3), the approximating normal distribution is not $N(0,1)$ but $N(0, e_n^2)$. Therefore, setting $a_n^{-1} = [e_n^{-1} - 1]$, we further expand the terms $\psi_{s,n}^*(u)$ involving the factor $e_n$ to write

(3.4)
$$\psi_{s,n}^*(u) = \left[1 + \sum_{k=0}^{s-2}\sum_{j=0}^{\nu_k(s)}\sum_{i=0}^{r_{j,k}(s)} b_n^{-k/2} n^{-j/2} a_n^{-i} p_{i,j,k,n}(u)\right]\phi(u) + O(\delta_n)$$

$$\equiv \psi_{s,n}(u) + O(\delta_n) \quad \text{(say)}$$

uniformly in $u \in \mathbb{R}$ for some polynomials $p_{i,j,k,n}(\cdot)$'s with $p_{0,0,0,n}(\cdot) \equiv 0$, where $r_{j,k}(s) = \min\{i \in \mathbb{Z}_+ : b_n^{-k/2} n^{-j/2} a_n^{-(i+1)} = o(b_n^{-(s-2)/2})\}$. The $(s-2)$*th-order*

EE $\Psi_{s,n}$ for $T_n$ is defined as

$$\Psi_{s,n}(u_0) = \int_{-\infty}^{u_0} \psi_{s,n}(u)\,du, \qquad u_0 \in \mathbb{R}.$$

3.2. *Main results.* We are now ready to state the main result of this section.

THEOREM 3.1. *Suppose that assumptions* (A.1)–(A.6) *hold. Then*

$$\sup_{u \in \mathbb{R}} |P(T_n \leq u) - \Psi_{s,n}(u)| = O(b_n^{-(s-2)/2}(\log n)^{-2}) \qquad \text{as } n \to \infty.$$

Theorem 3.1 establishes validity of a general-order EE for the studentized statistic $T_n$. From (3.4), it is clear that the EE in the dependent case has contributions from three sources, namely a series in powers of $n^{-1/2}$ coming from the estimation of $\theta$ (by $\hat{\theta}_n$), a series in powers of $b_n^{-1/2}$, coming from the stochastic error in estimating the asymptotic variance $\tau_{1n}^2$ of $\hat{\theta}_n$, and a third series in powers of $a_n^{-1}$, coming from the bias of the estimator $\hat{\tau}_n^2$. This should be compared with the relatively simpler form of the EE in the i.i.d. case [cf. (1.7)], where one is concerned with the estimation of *finite-dimensional* parameters only. In general, it is not possible to simplify the EE any further. Since the form of the EE is rather complicated, to provide further insight, we now provide an explicit expression for the *third-order EE* $\Psi_{s,n}(\cdot)$ with an error rate $o(n^{-1})$ in Section 3.3 below.

REMARK 3.1. The arguments developed here can be easily adapted to establish valid EEs for studentized versions of a slightly more general class of estimators than that given by the smooth function model [cf. (1.3)]. Specifically, let $\tilde{\theta}_n$ be an estimator of $\theta$ such that

$$(3.5) \qquad n^{1/2}(\tilde{\theta}_n - \theta) = \gamma' \bar{X}_n + \sum_{k=1}^{a(s)} n^{-k/2} \tilde{p}_{kn}(\bar{X}_n) + R_n$$

for some (nonrandom) $\gamma \in \mathbb{R}^d \setminus \{0\}$ and for some polynomials $\tilde{p}_{kn}(\cdot)$ with bounded coefficients, where $R_n$ satisfies the negligibility condition:

$$(3.6) \qquad P(|R_n| > Cb_n^{-(s-2-\kappa)}) = O(\delta_n)$$

for some $\kappa \in (0,1)$. Under suitable smoothness conditions on the score function $\psi(\cdot)$, the M-estimators in (1.4) admit such a representation, and hence the EE techniques developed here can be readily applied to studentized M-estimators in (1.4).

3.3. *Explicit expression for the third-order EE.* To derive the explicit EE formula with an error of order $o(n^{-1})$, note that under (A.1), the terms of order $b_n^{-r/2}$, $r \geq 4$ are $o(n^{-1})$. Hence, we restrict attention to the $(s-2)$th-order EE for $T_n$ with $s = 5$. For $s = 5$, the constants $\nu_k(s)$ in (3.1) are given by $\nu(s) = \nu_0(s) = 2$, $\nu_1(s) = 1$ and $\nu_k(s) = 0$ for $k = 2, 3$. The explicit forms of the EEs for $T_{1n}$ can now be obtained by carrying out the steps (3.1)–(3.4). These steps involve long and tedious algebraic computations and are given in in Section 6.3. Here, we state the resulting formulae only. To that end, let $H_k(\cdot)$ denote the $k$th-order Hermite polynomial, $k \geq 0$, defined by

$$H_k(x)\phi(x) = (-1)^k \frac{d^k}{dx^k} \phi(x), \qquad x \in \mathbb{R}.$$

PROPOSITION 3.2. *Suppose that assumptions* (A.1)–(A.6) *of Section 2 hold with* $s = 5$ *if* $\ell = O(n^{1/3})$ *and with* $s = 6$ *if* $\ell = O(n^{1/2-\kappa})$. *Then, the density of the preliminary EE* $\psi_n^*$ *[cf.* (3.3)*] for* $T_n$ *with an error of order* $o(n^{-1})$ *is given by*

$$\psi_n^*(x) = \phi_{e_n}(x)[1 + n^{-1/2}q_{1n}(x) + b_n^{-1}q_{2n}(x) \tag{3.7}$$
$$+ n^{-1}q_{3n}(x) + b_n^{-3/2}q_{4n}(x)], \qquad x \in \mathbb{R},$$

*where the polynomials* $q_{in}$'s *are as in* (6.57). *If, in addition,* $a_n^{-1} = O(n^{-1/3})$, *then*

$$\sup_{x \in \mathbb{R}} |P(T_n \leq x) - \Psi_{n,s}(x)| = O(n^{-1}(\log n)^{-2}),$$

*where, for* $x \in \mathbb{R}$,

$$\Psi_{n,s}(x) = \Phi(x) + \left[\left\{\frac{x}{a_n} - \frac{x^3}{2a_n^2} + \frac{x^3 H_2(x)}{6a_n^3}\right\}\right.$$
$$+ \{n^{-1/2}p_{1n}(x) + n^{-1/2}a_n^{-1}p_{2n}(x)\}$$
$$\tag{3.8}$$
$$+ \{b_n^{-1}p_{3n}(x) + b_n^{-1}a_n^{-1}p_{4n}(x)\}$$
$$\left. + n^{-1}p_{5n}(x) + b_n^{-3/2}p_{6n}(x)\right]\phi(x),$$

*and where the coefficients of the polynomials* $p_{in}$'s *are* $O(1)$ *as* $n \to \infty$. *The exact expressions for* $p_{in}$'s *are given in* (6.58).

Now, we discuss some of the implications of Proposition 3.2. The form of the preliminary EE for $T_n$ shows that the terms of order $b_n^{-1/2}$ do not appear in the EE, as also are the terms of order $b_n^{-1/2}n^{-1/2}$. As a result, interactions between the $n^{-1/2}$ series and the $b_n^{-1/2}$-series show up only in the terms of order *smaller* than $b_n^{-1/2}n^{-1/2}$. As a consequence, existing results in the literature on second and third-order EEs for the studentized statistics fail to describe this interaction. Next, consider relation (3.8) which gives the explicit third-order EE for the studentized statistic $T_n$, assuming $a_n^{-1} = O(n^{-1/3})$. For most commonly used weights $w_{kn}$'s [cf. discussion of assumptions (A.1) and (A.2), Section 2], the bias of the studentizing factor $\hat{\tau}_n^2$ is $O(n^{-1/3})$, so that $a_n^{-1} = O(n^{-1/3})$, and hence, this additional condition is not very restrictive. It is clear from (3.8) that it is the $a_n^{-1}$ series, but *not* the $b_n^{-1/2}$ series, that may slow down the rate of Normal approximation to the distribution of $T_n$ from the standard rate $O(n^{-1/2})$ observed in the i.i.d. case. Consequently, large sample confidence intervals for the parameter $\theta$ based on $T_n$ can have *very poor coverage probability under dependence*, unless the weights $w_{kn}$'s are chosen carefully to make the bias small.

REMARK 3.2. The explicit EE results of Proposition 3.2 can also be used to determine the *optimal block size $\ell$* for minimizing the error of Normal approximation to the distribution of the studentized statistic $T_n$. Suppose that the weights $w_{kn}$'s in the definition of the studentizing factor $\hat{\tau}_n^2$ in (1.5) is generated by a bounded taper function $w:[0,1] \to \mathbb{R}$ as $w_{kn} = w(k/\ell)$, $1 \leq k \leq \ell$ [cf. (2.4)], with $w(0+) = 1$. Let $r$ be the *characteristic exponent* of $w(\cdot)$, i.e., $r$ be the largest integer such that

$$(3.9) \qquad w^{(r)} \equiv \lim_{y \to 0+} \frac{1 - w(y)}{y^r}$$

exists and is nonzero. Then, it can be shown that the bias of $\hat{\tau}_n^2$ is of the order of $\ell^{-r}$ and hence, $a_n^{-1} = c_r \ell^{-r}(1 + o(1))$ as $n \to \infty$, for some constant $c_r$ [involving the constant $w^{(r)}$ and derivatives of the spectral density function of the process $\{Z_i^*\}_{i \in \mathbb{Z}}$ where $Z_i^* = h(\mu)'X_i$]. From (3.8), it now follows that the rate of Normal approximation is optimized provided $a_n^{-1}$ decays like $n^{-1/2}$, which in turn implies that the *optimal* choice of the block size is of the order $n^{1/(2r)}$. Thus, for the Bartlett lag-window, $w(y) = (1 - y)$, $y \in [0, 1]$, we have $r = 1$ in (3.9) and the optimal block size is of the order $n^{1/2}$. Similarly, for the studentization factor given by the Moving Block Bootstrap (MBB) method (cf. Section 4.2 below), the weights are asymptotically equivalent to those generated by the Bartlett lag-window and hence, the optimal choice of the block size here is also of the order $n^{1/2}$. [Although assumption (A.1) rules out blocks of size $n^{1/2}$ for a general-order EE, a second-order EE can be established for $\ell = O(n^{1/2})$ with some additional work, which shows that these are indeed the optimal order for the case $r = 1$.] Next, consider the Daniell window, $w(y) = \frac{\sin \pi y}{\pi y}$, and the Parzen window, $w(y) = [1 - 6y^2 + 6y^3]\mathbb{1}(0 \leq y \leq 1/2) + 2(1 - y)^3\mathbb{1}(1/2 \leq y \leq 1)$. For both of these, $r = 2$ and hence, the optimal block size is of the order $n^{1/4}$. Finally, consider the flat-top kernels of Politis and Romano (1995). For these, $w^{(r)} = 0$ for all $r \geq 1$, and we have what are called the *infinite-order* windows. In such cases, under (A.1)–(A.5), the bias of $\hat{\tau}_n^2$ decays exponentially fast, and the optimal order of $\ell$ is $O(\log n)$.

**4. EEs under alternative studentizations.** In this section, we give results on general-order EEs for studentized statistics for various alternative choices of the studentizing factor that have been proposed in the literature.

4.1. *Lag covariance estimators with a different scaling.* First, consider the effects of replacing $\hat{\Gamma}_n(k)$ by

$$(4.1) \quad \hat{\Gamma}_n^{[1]}(k) = (n-k)^{-1} \sum_{i=1}^{n-k} (X_i - \bar{X}_n)(X_{i+k} - \bar{X}_n)', \qquad 0 \leq k \leq \ell,$$

i.e., we replace the scaling factor $n^{-1}$ in $\hat{\Gamma}_n(k)$ by $(n-k)^{-1}$. Let

$$T_n^{[1]} = \sqrt{n}(\hat{\theta}_n - \theta)/\hat{\tau}_n^{[1]},$$

where $\hat{\tau}_n^{[1]}$ is defined by replacing $\hat{\Gamma}_n(k)$'s in the definition of $\hat{\tau}_n$ by $\hat{\Gamma}_n^{[1]}(k)$, $0 \leq k \leq \ell$. In this case, the basic EE results of Section 3 readily yield EEs for $T_n^{[1]}$. Note that

$$\hat{\Gamma}_n^{[1]}(0) + \sum_{k=1}^{\ell} w_{kn}(\hat{\Gamma}_n^{[1]}(k) + \hat{\Gamma}_n^{[1]}(k)') = \hat{\Gamma}_n(0) + \sum_{k=1}^{\ell} \tilde{w}_{kn}(\hat{\Gamma}_n(k) + \hat{\Gamma}_n(k)'),$$

where $\tilde{w}_{kn} \equiv w_{kn}[n/(n-k)], 1 \leq k \leq \ell$. Therefore, it follows that an $(s-2)$th-order EE for $T_n^{[1]}$ *is valid under the same regularity conditions as those needed in the case of* $T_n$ with $\tilde{w}_{kn}$'s replacing the $w_{kn}$'s, and the corresponding EE is given by $\Psi_{s,n}(\cdot)$, again with $\tilde{w}_{kn}$'s replacing the $w_{kn}$'s.

4.2. *Block bootstrap variance estimators.* A popular approach to finding estimators of the standard error of $\hat{\theta}_n$ is through one of the many block bootstrap and subsampling methods proposed in the literature [cf. Lahiri (2003)]. In this section, we consider higher-order expansions for the studentized statistic, when the studentizing factor is generated by different block resampling methods. Let $\ell$ be an integer satisfying (A.1) and let $N = n - \ell + 1$. Here, $N$ denotes the number of overlapping blocks of size $\ell$ contained in $\{X_1, \ldots, X_n\}$. For $i = 1, \ldots, N$, define $V_i = \sum_{k=i}^{i+\ell-1} w_{kn}(X_k - \bar{X}_n)$, where $w_{kn}$'s are nonrandom weights satisfying $\sum_{k=1}^{\ell} w_{kn}^2 = 1$. Then, the *block resampling estimator* of $\text{Cov}(Z_n)$ (where $Z_n = n^{1/2}[\bar{X}_n - \mu]$) is given by

$$(4.2) \qquad \hat{\Sigma}_n^{\text{BR}} \equiv N^{-1} \sum_{i=1}^{N} V_i V_i'.$$

Note that for $w_{kn} \equiv \ell^{-1/2}$, this gives the MBB estimator of $\text{Cov}(Z_n)$ [Künsch (1989), Liu and Singh (1992)]. Due to the linearity of the sample mean in the $X_i$'s, $\hat{\Sigma}_n^{\text{BR}}$ with $w_{kn} \equiv \ell^{-1/2}$ is also the overlapping subsampling estimator of $\text{Cov}(Z_n)$ [cf. Politis and Romano (1994)]. Further, for suitable choices of $w_{kn}$'s, $\hat{\Sigma}_n^{\text{BR}}$ gives the tapered block bootstrap estimator of $\text{Cov}(Z_n)$ [cf. Paparoditis and Politis (2001)].

With the block resampling estimator $\hat{\Sigma}_n^{\text{BR}}$ of $\text{Cov}(Z_n)$ given by (4.2), we now define the corresponding studentized version of $\hat{\theta}_n$ as

$$T_n^{[2]} = \sqrt{n}(\hat{\theta}_n - \theta)/\hat{\tau}_n^{[2]},$$

where $\hat{\tau}_n^{[2]} = \max\{n^{-1}, h(\bar{X}_n)'[\hat{\Sigma}_n^{\text{BR}}]h(\bar{X}_n)\}^{1/2}$. A general-order EE for $T_n^{[2]}$ continue to hold under regularity conditions similar to (A.1)–(A.6). To state the result, define

$$W_{jn}^{[2]} = \left(n^{1/2}[\bar{X}_n - \mu]'; \left[\text{VEC}\left(\ell^{-1}\sum_{i=(j-1)\ell+1}^{j\ell \wedge n} U_{1i}U_{1i}\right)\right]'\right)',$$

(4.3)
$$j = 1, \ldots, b_{0n},$$

where $U_{1i} = nN^{-1}\sum_{k=i}^{i+\ell-1} w_{(k-i+1)n}(X_k - \mu)$, $i = 1, \ldots, N$ and $U_{1i} = 0$ for $i = N+1, \ldots, n$. Here and in the following, for a $d \times d$ matrix $A = ((a_{ij}))$, $\text{VEC}(A)$ denotes the $d_1 = d(d+1)/2$-dimensional vector, given by $\text{VEC}(A) = (a_{11}, \ldots, a_{1d}; a_{22}, \ldots, a_{2d}; \ldots; a_{dd})'$. The collection of random vectors $\{W_{jn}^{[2]} : j = 1, \ldots, b_{0n}\}$ gives the basic block variables for the studentized statistic $T_n^{[2]}$ with the block resampling studentizing factor $\hat{\tau}_n^{[2]}$. The following result gives the conditions for a valid $(s-2)$th-order EE for $T_n^{[2]}$.

THEOREM 4.1. *Suppose that assumption* (A.1)–(A.6) *hold with $W_{jn}$'s replaced by $W_{jn}^{[2]}$'s. Then*

$$\sup_{u \in \mathbb{R}}|P(T_n^{[2]} \leq u) - \Psi_{s,n}^{[2]}(u)| = O(\delta_n) \qquad as\ n \to \infty,$$

*where $\Psi_{s,n}^{[2]}$ is obtained from the definition of $\Psi_{s,n}$ by replacing $\hat{A}_{kn}, k = 1, 2, 3$ with $\hat{A}_{kn}^{[2]}, k = 1, 2, 3$ in (3.1) and where $\hat{A}_{kn}^{[2]}, k = 1, 2, 3$ are as defined in (6.59). Further, for the linear process in (2.5), assumptions* (A.3)–(A.6) *[with $W_{jn} = W_{jn}^{[2]}$ in* (A.6)*] continue to hold solely under the conditions of Proposition 2.1.*

4.3. *Variance estimators without mean correction.* In some limited applications, one can assume that the mean $\mu$ of the $X_1$'s is *known*. For example, if the $X_i$'s are univariate and one is interested in testing hypotheses of the form $H_0: \mu = \mu_0$ against $H_1: \mu \neq \mu_0$ for some given $\mu_0$, it may be reasonable to use a studentizing factor without *explicit estimation* of the mean parameter $\mu$. In such situations, one would simply make use of the hypothesized value $\mu_0$ of $\mu$ to define a test-statistic. This often simplifies the proof and derivation of the EE for the studentized statistic.

For the rest of this subsection, we suppose that $\mu$ is *known*. Further, for notational simplicity, we set $\mu = 0$ (otherwise, replace $X_i$ in each occurrence with $X_i - \mu$). Then, each of the studentizing factors described above leads to a simpler form where the $\bar{X}_n$'s are dropped from the corresponding formulae for $\hat{\tau}_n^{[r]}$'s. As a specific case, we consider the known-mean version of the studentizing factor arising from $\hat{\tau}_n^{[2]}$, given by

$$\hat{\tau}_n^{[3]} = \max\left\{n^{-1}, h(\bar{X}_n)'\left[N^{-1}\sum_{i=1}^{N} U_{1i}U_{1i}'\right]h(\bar{X}_n)\right\}^{1/2},$$

where, from above, $U_{1i} = nN^{-1}\sum_{k=i}^{i+\ell-1} w_{(k-i+1)n}X_k$, $i = 1, \ldots, N$ and $U_{1i} = 0$ for $i = N+1, \ldots, n$. Set $T_n^{[3]} = \sqrt{n}(\hat{\theta}_n - \theta)/\hat{\tau}_n^{[3]}$. Then we have the following result.

THEOREM 4.2. *Suppose that assumption* (A.1)–(A.5) *hold with $W_{jn}$'s replaced by $W_{jn}^{[2]}$ [cf. (4.3)], and the following weaker version of* (A.6) *holds:*

(A.6)′ *There exists a constant $\kappa \in (0, 1)$ and a sequence $\{d_n\} \subset [1, \infty)$ with $d_n = O(\ell)$ and $d_n^2 = O(b_n^{1-\kappa})$ such that for all $n \geq \kappa^{-1}$ and $\kappa^{-1} < m \leq [\log n]^{1/\kappa}$,*

$$(4.4) \quad \max_{j_0 \in J_n} \sup_{t \in B_n} E|E\{\exp(\iota t'[W_{(j_0-m)n}^{[2]} + \cdots + W_{(j_0+m)n}^{[2]}])|\tilde{\mathcal{D}}_{j_0}\}| \leq 1 - \kappa,$$

*where $J_n$, $B_n$, and $\tilde{\mathcal{D}}_{j_0}$ are as in* (A.6). *Then*

$$\sup_{u \in \mathbb{R}}|P(T_n^{[3]} \leq u) - \Psi_{s,n}^{[3]}(u)| = O(\delta_n) \qquad \text{as } n \to \infty,$$

*where $\Psi_{s,n}^{[3]}$ is obtained from the definition of $\Psi_{s,n}$ by replacing $\hat{A}_{kn}, k = 1, 2, 3$ with $\hat{A}_{kn}^{[3]}, k = 1, 2, 3$ in (3.1) and where $\hat{A}_{1n}^{[3]} = \hat{A}_{1n}^{[2]}$, with $\hat{A}_{1n}^{[2]}$ as in (6.59), and where $\hat{A}_{kn}^{[3]} = 0$ for $k = 2, 3$. Further, for the linear process in (2.5), assumptions* (A.3)–(A.5) *and* (A.6)′ *hold solely under the conditions of Proposition 2.1.*

Note that Theorem 4.2 gives a valid EE for $T_n^{[3]}$ under a slightly weaker conditional Cramer condition (where $E\sup_{t \in B_n}|\cdot|$ is replaced by $\sup_{t \in B_n} E|\cdot$

|). This is a consequence of not using an explicit mean correction in the studentizing factor $\hat{\tau}_n^{[3]}$. In this case, valid EE can be derived by combining the EE results of [L] for sums of block variables and the transformation technique of Bhattacharya and Ghosh (1978). No iterated conditioning arguments, like those used in the proof of Theorem 3.1 (cf. Section 5 below) are needed. Velasco and Robinson (2001) proved a third-order EE for a mean-uncorrected version of the studentized statistic, for Gaussian processes. Theorem 4.2 supplements their results by dropping the Gaussianity requirement and giving a general-order EE, under regularity conditions (A.1)–(A.5) and (A.6)′.

**5. Edgeworth expansions for perturbed sums of block variables.** It can be shown (cf. Section 6) that the basic building block in the stochastic approximation to the studentized statistics $T_n$ is given by a perturbed sum of block variables of the form:

(5.1) $$R_n^* \equiv \left( n^{1/2} \bar{X}_n', \frac{1}{\sqrt{n\ell}} \sum_{i=1}^{b_n} Y_{in}' + b_n^{-1} \xi_n' \right)',$$

where, $Y_{in}$'s are as in Section 2 and $\xi_n$ is a quadratic form in $n^{1/2}\bar{X}_n$ and the boundary blocks of length $\ell$ [cf. (6.2) below]. As a result, for deriving the EEs in Section 4, we need to develop an EE theory for $R_n^*$, which we address in this section.

Note that the second term in $R_n^*$ has a nonstandard form which gives rise to some difficult technical issues that must be resolved for deriving valid EEs for such a statistic. Specifically, note that with $d = 1$ (for notational convenience), $\xi_n$ is a variable of the form

$$\xi_n = n^{1/2} \bar{X}_n [_3W_{1n} + {}_3W_{b_{0n}n}],$$

where $_3W_{jn}$, $j = 1, b_{0n}$ are certain boundary block variables (based on blocks of size $\ell$) that are asymptotically normal. Although it may appear reasonable at the first glance, deriving an EE for $R_n^*$ from the joint distribution of $n^{1/2}\bar{X}_n$, $\frac{1}{\sqrt{n\ell}} \sum_{i=1}^{b_n} Y_{in}$ and $_3W_{jn}$, $j = 1, b_{0n}$ is not very useful for our purpose, as the EEs for $_3W_{jn}$, $j = 1, b_{0n}$ are in powers of $\ell^{-1/2}$, which can be too slow to give an overall error rate of $o(b_n^{-(s-2)/2})$, particularly when $\ell \sim (\log n)^C$. This observation leads us to consider finding an EE for $\xi_n$ directly. But by taking the direct approach, we have to contend with a complication that arises from the fact that the boundary block variables, being *common multipliers* of all the $X_i$'s in $n^{1/2}\bar{X}_n$, induce strong dependence of among the product variables. To deal with this, here we develop a conditioning argument with resect to suitable initial and end segments of the variables $X_1, \ldots, X_n$ and derive a valid EE for the joint distribution of $n^{1/2}\bar{X}_n$ and $\frac{1}{\sqrt{n\ell}} \sum_{i=1}^{b_n} Y_{in} + b_n^{-1}\xi_n$. A representation result of Bradley (1983) (cf. Lemma 6.1 below) for strongly mixing random variables plays a crucial role in carrying out the conditioning step in the proof.

Next, we define the $(s-2)$th-order EE $\Upsilon_{s,n}(\cdot)$ for $R_n^*$ through its Fourier transform $\hat{\Upsilon}_{s,n}(t) \equiv \int \exp(\iota t'x)\Upsilon(dx)$, $t \in \mathbb{R}^{d_0}$. Let $S_n^{[0]}(t) = t'R_n^*$, $t \in \mathbb{R}^{d_0}$. Note that $ES_n^{[0]}(t) = b_n^{-3/2} t' E\xi_{1n}$ and $\text{Var}(S_n^{[0]}(t)) = t'\Xi_n t + 2b_n^{-3/2} t' ES_n \xi_{1n}' t + b_n^{-3} t' E\xi_{1n}\xi_{1n}' t$, where $\xi_{1n} = (0':\xi_n')'$. Also, by Lemma 6.3, for $3 \leq r \leq s$, the $r$th-order cumulant [cf. (6.4)] of $S_n^{[0]}(t)$ is

$$\mathcal{K}_0^{[0]}([S_n^{[0]}(t)]^{\diamond r}) = \sum_{|a_1|=1} \cdots \sum_{|a_r|=1} t^{a_1+\cdots+a_r} \mathcal{K}_0^{[0]}(S_n^{[0]}(a_1),\ldots,S_n^{[0]}(a_r))$$

$$= b_n^{-(r-2)/2} \sum_{|\alpha|=r} t^\alpha \left( \sum_{j=0}^{s-r} b_n^{-j/2} \lambda_{j,n}^{(\alpha)} \right)$$

for some constants $\lambda_{j,n}^{(\alpha)} \in \mathbb{R}$, with $\lambda_{j,n}^{(\alpha)} = O(1)$ for all $j$ and $\alpha$. Next, define the EE polynomials $\{q_{r,n}x : r = 1,\ldots, s-2\}$ by the identity:

$$\exp\Bigg( b_n^{-3/2}[\iota t' E\xi_{1n} - t' ES_n\xi_{1n}' t] - \frac{b_n^{-3}}{2} t' E\xi_n\xi_{1n}' t$$
$$+ \sum_{r=3}^{s} \frac{b_n^{-(r-2)/2}}{r!} \sum_{|\alpha|=r} (\iota t)^\alpha \left[ \sum_{j=0}^{s-r} b_n^{-j/2} \lambda_{j,n}^{(\alpha)} \right] \Bigg)$$
$$= 1 + \sum_{r=1}^{s-2} b_n^{-r/2} q_{r,n}(\iota t) + o(b_n^{(s-2)/2}).$$

The Fourier transform $\hat{\Upsilon}_{s,n}$ of $\Upsilon_{s,n}(\cdot)$ is defined as

$$\hat{\Upsilon}_{s,n}(t) = \exp(-t'\Xi_n t/2)\left[1 + \sum_{r=1}^{s-2} b_n^{-r/2} q_{r,n}(\iota t)\right], \qquad t \in \mathbb{R}^{d_0}.$$

With this, we now have the main result of this section.

THEOREM 5.1. *Let $R_n^*$ be as defined in (5.1). Then under assumptions* (A.1)–(A.6),

(5.2) $$\sup_{B \in \mathcal{B}} |P(R_n^* \in B) - \Upsilon_{s,n}(B)| = O(b_n^{-(s-2)/2}[\log n]^{-2}),$$

*where $\mathcal{B}$ is a class of Borel subsets of $\mathbb{R}^{d_0}$ satisfying*

$$\sup_{B \in \mathcal{B}} \Phi_\Xi([\partial B]^\varepsilon) = O(\varepsilon) \qquad \text{as } \varepsilon \downarrow 0.$$

**6. Proofs.** For clarity of exposition, we shall switch the order of proofs of the results from Sections 2, 3 and 4 and those from Section 5. In Section 6.1, we set forth some common notation and state some relevant facts that are used in the rest of the section. The proof of Theorem 5.1 for a general $d \geq 1$ and various intermediate steps (cf. Lemmas 6.1–6.6) are presented in Section 6.2. Proofs of the results from Sections 2, 3 and 4 are given in Sections 6.3–6.5, respectively.

6.1. *Notation and preliminaries.* With a given $d \geq 1$, let $d_0 = d_1 + d$, where $d_1 = d(d+1)/2$. Thus, $d_0$ is the dimension of the block variables $W_{jn}$'s in (2.2). For a $d \times d$ matrix $A = ((a_{i,j}))$, let $\text{VEC}(A)$ be the $d_1$-dimensional vector, consisting of the on- and above-the-diagonal entries of $A$ (with the lexicographic ordering). Let $\text{SVEC}(A) = \text{VEC}(A) + \text{VEC}(A')$. Next, define the reverse operation $\text{VEC}^{-1}(\cdot)$ from $\mathbb{R}^{\Lambda_0}$ into the class of $d \times d$ matrices by setting

$$[\text{VEC}^{-1}(x)]_{i,j} = x_{i,j} \in ((i,j) \in \Lambda_0).$$

Thus, $\text{VEC}^{-1}(x)$ is an upper triangular matrix with the elements from $x$ filling in the on- and above-the-diagonal positions and with zeros below the diagonal. Also, set $\text{SVEC}^{-1}(x) = \text{VEC}^{-1}(x) + [\text{VEC}^{-1}(x)]'$. Note that for $t \in \mathbb{R}^{\Lambda_0}$ and $x, y \in \mathbb{R}^d$,

$$
\begin{aligned}
t' \text{VEC}(xy') &= x'[\text{VEC}^{-1}(t)]y, \\
(6.1) \qquad t' \text{SVEC}(xy') &= x'[\text{SVEC}^{-1}(t)]y \quad \text{and} \\
\|\text{SVEC}^{-1}(t)x\| &\leq C(d)\|t\|\|x\|.
\end{aligned}
$$

Here and in the proofs below, let $C, C(\cdot)$ denote generic constants in $(0, \infty)$, depending on their arguments (if any), but not on $n$ and $\ell$. For notational simplicity, we will suppress the dependence of $C, C(\cdot)$ on the dimension variables $d_0, d_1$ and $d$ and the constant $\kappa$ appearing in assumptions (A.1)–(A.6). Let $e_1, \ldots, e_d$ denote the standard basis of $\mathbb{R}^d$. Let $\psi : [0, \infty) \to [0, \infty)$ be an infinitely differentiable nondecreasing function satisfying $\psi(u) = u$ for $0 \leq u \leq 1$ and $\psi(u) = 2$ for $u \geq 2$. For $c > 0$, define the truncation function $g(\cdot; c) : \mathbb{R}^k \to \mathbb{R}^k$ (for a $k \in \mathbb{N}$) by

$$g(x; c) = \begin{cases} \dfrac{cx}{\|x\|} \cdot \psi\left(\dfrac{\|x\|}{c}\right), & \text{if } x \neq 0, \\ 0, & \text{if } x = 0. \end{cases}$$

Thus, $g(x; c) = x$ for $\|x\| \leq c$ and $\|g(x; c)\| \leq 2c$ for all $x \in \mathbb{R}^k$. For $n \geq 1$, set $c_n = b_n^{1/2}(\log(n+1))^{-2}$ and let $\delta_{n,C} = b_n^{-(s-2)/2}(\log[n+1])^{-C}$ for some $C \in (0, \infty)$. Note that $\delta_n = \delta_{n,2}(1 + o(1))$. For a sequence of random variables $\{W_n\}_{n\geq 1}$ and a sequence real numbers $\{x_n\}_{n\geq 1} \in (0,\infty)$, we write $W_n = \tilde{O}_p(x_n)$ if

$$P(|W_n| > Cx_n) = O(\delta_n) \quad \text{for some } C > 0.$$

For notational simplicity, we will often drop the subscript $n$ and write $b_n = b$, $b_{0n} = b_0, \ell_n = \ell$, etc., whenever this convention does not create confusion. Also, unless otherwise stated, we set $\mu = 0$ in all of the steps below. Thus, the $X$-variables appearing in the expressions below are centered; in particular, $\bar{X}_n = O_p(n^{-1/2})$.

Let $W_{jn} = ({}_1W'_{jn}, {}_2W'_{jn})'$, where ${}_1W_{jn}$ is $d \times 1$. Similarly, partition $S_n \equiv b_n^{-1/2} \sum_{j=1}^{b_{0n}} W_{jn}$ as $S_n = ({}_1S'_n, {}_2S'_n)'$. Let

(6.2) $$\xi_n \equiv n^{1/2} \operatorname{SVEC}(\bar{X}_n[{}_3W_{1n} + {}_3W_{b_0n}]'),$$

where ${}_3W_{1n} = \ell^{-3/2} \sum_{i=1}^{\ell} (\sum_{k=i}^{\ell} w_{kn}) X_i$ and ${}_3W_{b_0n} = \ell^{-3/2} \sum_{i=1}^{\ell} (\sum_{k=i}^{\ell} w_{kn}) \times X_{n-i+1}$. Also, with ${}_3W_{jn} = 0$ for $j = 2, \ldots, b_{0n} - 1$, set

$$W_{jn}^{[0]} = ({}_1W'_{jn}; {}_2W'_{jn}; {}_3W'_{jn})', \quad j = 1, \ldots, b_{0n}.$$

Unless otherwise indicated, a check (ˇ) and a tilde (˜) on a random vector (at stage $n$ in the asymptotic argument) denote its truncated version (at $c_n$) and its centered truncated version, respectively. Thus, ${}_r\check{W}_{jn} = g({}_rW_{jn}; c_n)$ and ${}_r\tilde{W}_{jn} = {}_r\check{W}_{jn} - E{}_r\check{W}_{jn}$ for $r = 1, 2, 3$ and $j = 1, \ldots, b_0$. Next (with a slight abuse of the above notation), define $\check{W}_{jn}^{[0]}$ and $\tilde{W}_{jn}^{[0]}$ by replacing $\{{}_rW_{jn} : r = 1, 2, 3\}$ with $\{{}_r\check{W}_{jn} : r = 1, 2, 3\}$ and $\{{}_r\tilde{W}_{jn} : r = 1, 2, 3\}$, respectively. Set $\tilde{\xi}_n = \operatorname{SVEC}({}_1\tilde{S}_n[{}_3\tilde{W}_{1n} + {}_3\tilde{W}_{b_0n n}]')$ and $\tilde{S}_n^{[0]}(t) = t'\tilde{S}_n + b_n^{-3/2} t'\tilde{\xi}_{1n}$, with $\tilde{\xi}_{1n} = (0', \tilde{\xi}'_n)' \in \mathbb{R}^{d_0}$. Note that for $t = (t'_1, t'_2)' \in \mathbb{R}^d \times \mathbb{R}^{\Lambda_0}$,

(6.3)
$$\tilde{S}_n^{[0]}(t) = b_n^{-1/2} \sum_{j=1}^{b_{0n}} [t'\tilde{W}_{jn} + b_n^{-3/2}({}_3\tilde{W}_{1n} + {}_3\tilde{W}_{b_0n n})'[\operatorname{SVEC}^{-1}(t_2)]({}_1\tilde{W}_{jn})]$$

$$\equiv b_n^{-1/2} \sum_{j=1}^{b_{0n}} \tilde{V}_{jn}(t), \quad \text{say.}$$

Next, for an integrable random variable $V$ on $(\Omega, \mathcal{F}, P)$, define

$$E_t V = [H_n^{[0]}(t)]^{-1} E[V \exp(\tilde{S}_n^{[0]}(\iota t))], \quad t \in \mathbb{R}^{d_0},$$

where $H_n^{[0]}(t) = E \exp(\iota[\tilde{S}_n^{[0]}(t)])$, $t \in \mathbb{R}^{d_0}$. Also, let $H_n(t) = E \exp(\iota t'\tilde{S}_n)$, $t \in \mathbb{R}^d$. Define the semi-invariants $\mathcal{K}_t^{[0]}(V_1, \ldots, V_p)$ in variables $V_1, \ldots, V_p$ on $(\Omega, \mathcal{F}, P)$, by

$$\iota^p \mathcal{K}_t^{[0]}(V_1, \ldots, V_p)$$

$$
(6.4) \quad = \frac{\partial}{\partial \varepsilon_1} \cdots \frac{\partial}{\partial \varepsilon_p}\bigg|_{\varepsilon_1=\cdots=\varepsilon_p=0}
$$
$$
\times \log E \exp(\iota t' \tilde{S}_n^{[0]} + \iota[\varepsilon_1 V_1 + \cdots + \varepsilon_p V_p]), \qquad t \in \mathbb{R}^{d_0}.
$$

Then, it is a well-known fact that (cf. (3.13) of [GH])

$$
(6.5) \quad \mathcal{K}_t^{[0]}(V_1,\ldots,V_p) = \sum_{j=1}^p {\sum}^{*} c(I_1,\ldots,I_j) \prod_{i=1}^j E_t \prod_{k \in I_i} V_k,
$$

where $\sum^*$ extends over all disjoint unions $I_1 \cup \cdots \cup I_j = \{1,\ldots,p\}$ and where $c(I_1,\ldots,I_j)$'s are combinatorial coefficients. Also, for any nonempty and disjoint decomposition of $\{1,\ldots,p\}$ into the sets $J_1$ and $J_2$, one has (cf. (3.14) of [GH]),

$$
(6.6) \quad 0 = \sum_{j=1}^p {\sum}^{*} c(I_1,\ldots,I_j) \prod_{i=1}^j \left( E_t \prod_{k \in I_{i1}} V_k \right) \left( E_t \prod_{k \in I_{i2}} V_k \right),
$$

where $I_{ir} = I_i \cap J_r, 1 \leq i \leq j, r = 1,2$.

Next, write $\mathcal{K}_t^{[0]}(V_1^{\diamond p}, V_2^{\diamond q}) = \mathcal{K}_t^{[0]}(V_1,\ldots,V_1,V_2,\ldots,V_2)$ where, on the right-hand side, $V_1$ appears $p$-times and $V_2$ appears $q$ times. Expanding $\log H_n^{[0]}(t)$ in the Taylor's series, one has

$$
(6.7) \quad \log H_n^{[0]}(t) = \sum_{r=2}^s \iota^r \mathcal{K}_0^{[0]}([\tilde{S}_n^{[0]}(t)]^{\diamond r}) + R_n(t),
$$

where $R_n(t) = [\int_0^1 (1-\eta)^s \iota^{s+1} \mathcal{K}_{\eta t}^{[0]}([\tilde{S}_n^{[0]}(t)]^{\diamond(s+1)})\,d\eta]/s!$.

6.2. *Proof of the main result from Section 5.* Before proving Theorem 5.1, we state and prove some of the intermediate steps in its proof as lemmas. The first lemma is a variant of Bradley (1983)'s result on representation of strongly mixing random variables that plays a crucial role in the rest of the proof.

LEMMA 6.1. *Let $X$ be a $r$-random vector and let $Y = (Y_1,\ldots,Y_k)'$ be a $k$-random vector on a probability space with $\|Y\|_\gamma < \infty$ and with $\alpha_0 = \alpha(\sigma\langle X\rangle, \sigma\langle Y\rangle)$, where $\gamma \in (0,\infty)$ and $k,r \in \mathbb{N}$. Then there exist a $k$-random vector $Y^*$ (on a possibly enlarged probability space) such that (i) $Y$ and $Y^*$ have the same probability distribution on $\mathbb{R}^k$, (ii) $X$ and $Y^*$ are independent and (iii) for all $t \in (0, \min_{j=1,\ldots,k} \|Y_j\|_\gamma]$,*

$$
(6.8) \quad P\left( \max_{j=1,\ldots,k} |Y_j - Y_j^*| > t \right) \leq 18k \left( \alpha_0^2 \max_{j=1,\ldots,k} \|Y_j\|_\gamma \cdot t^{-1} \right)^{\gamma/(2+\gamma)}.
$$

PROOF. For $k = 1$, this is a restatement of Theorem 3 of Bradley (1983). For $k \geq 2$, one can use finite subadditivity and apply Bradley (1983)'s result term-wise to get (6.8). $\square$

LEMMA 6.2. *Let $Z_1, \ldots, Z_m$ be random vectors (of dimensions $k_1, \ldots, k_m \in \mathbb{N}$, respectively) on a probability space with $E\|Z_i\|^\gamma < \infty$ for some $\gamma \in (0, \infty)$. Let $\alpha_{0j} = \alpha(\sigma\langle Z_j \rangle, \sigma\langle Z_{j+1}, \ldots, Z_m \rangle)$ for $j = 1, \ldots, m-1$. Then there exist independent random vectors $\{Z_1^*, \ldots, Z_m^*\}$ (on a possibly enlarged probability space) such that* (i) *$Z_j$ and $Z_j^*$ have the same probability distribution on $\mathbb{R}^{k_j}$, $j = 1, \ldots, m$,* (ii) *$Z_m = Z_m^*$, and* (iii) *for all $j = 1, \ldots, m-1$, and $t \in (0, t_{0j}]$,*

$$(6.9) \quad P\Big(\max_{i=1,\ldots,k_j} |Z_{ji} - Z_{ji}^*| > t\Big) \leq 18k_j \Big(\alpha_{0j}^2 \max_{i=1,\ldots,k_j} \|Z_{ji}\|_\gamma \cdot t^{-1}\Big)^{\gamma/(2+\gamma)},$$

*where $t_{0j} = \min\{\|Z_{ji}\|_\gamma : i = 1, \ldots, k_j\}$, and $Z_{ji}, i = 1, \ldots, k_j$ are the components of $Z_j$.*

PROOF. First, take $Y = Z_1$ and $X = (Z_2', \ldots, Z_m')'$ in Lemma 6.1 to get a copy $Z_1^*$ of $Z_1$ that is independent of $\{Z_2, \ldots, Z_m\}$ and that satisfies (6.9) with $j = 1$. Next, set $Y = Z_2$ and $X = (Z_3', \ldots, Z_m'; Z_1^{*\prime})'$ in Lemma 6.1. Using the independence of $Z_1^*$ and $Z_2, \ldots, Z_m$, it is easy to check that

$$\alpha(\sigma\langle Y \rangle, \sigma\langle X \rangle) = \sup_{A,B} |P(Z_2 \in A, X \in B) - P(Z_2 \in A)P(X \in B)|$$

$$\leq E\Big(\sup_{A,B} |P_{\cdot|Z_1^*=z_1}(Z_2 \in A, (Z_3', \ldots, Z_m'; z_1)' \in B)$$

$$- P(Z_2 \in A)P_{\cdot|Z_1^*=z_1}((Z_3', \ldots, Z_m'; z_1)' \in B)|\Big)$$

$$\leq \alpha_{0j}.$$

Hence, by Lemma 6.1, there exists a copy $Z_2^*$ of $Z_2$ that is independent of $X = (Z_3', \ldots, Z_m'; Z_1^{*\prime})'$ and that satisfies (6.9) with $j = 2$. Now repeat this process $(m-1)$-times to get independent variables $Z_1^*, \ldots, Z_{m-1}^*$ satisfying (6.9) that are *all* independent of $Z_m$. Finally, set $Z_m^* = Z_m$ to complete the proof of the lemma. $\square$

REMARK 6.1. Note that in Lemma 6.2, the bound in (6.9) involves the dimensions of the *first $(m-1)$* vectors only, but *not* the dimension of $Z_m$.

LEMMA 6.3. *Suppose that assumptions* (A.1)(i) *and* (A.2)–(A.5) *hold. Then, for any $2 \leq r \leq s$ $a_1, \ldots, a_r \in \mathbb{R}^{d_0}$ with $\|a_j\| \leq 1$, $j = 1, \ldots, r$,*

$$(6.10) \quad \mathcal{K}_0(\tilde{S}_n^{[0]}(a_1), \ldots, \tilde{S}_n^{[0]}(a_r)) \leq C b_n^{-(r-2)/2},$$

*uniformly in $a_1, \ldots, a_r \in \mathbb{R}^{d_0}$ with $\|a_j\| \leq 1$, $j = 1, \ldots, r$.*

PROOF. Fix $2 \leq r \leq s$ and $a_1, \ldots, a_r \in \mathbb{R}^{d_0}$ with $\|a_j\| \leq 1$. Let $s_r = \min\{\lfloor (s-2)/3 \rfloor, r\}$. Note that for $r \geq \lfloor (s-2)/3 \rfloor$, $b_n^{-3(s_r+1)/2} = o(b_n^{-[s-2]/2})$. Also, note that for any random variables $V_1, \ldots, V_r$,

$$|\mathcal{K}_0(V_1, \ldots, V_r)| \leq C(r) \prod_{j=1}^{r} \{E\|V_j\|^r\}^{1/r} \quad \text{and} \quad E\|\tilde{S}_n\|^r + E\|\tilde{\xi}_{1n}\|^r = O(1).$$

Hence, using the multilinearity property of the semiinvariants, we have

$$\begin{aligned}
\mathcal{K}_0(\tilde{S}_n^{[0]}(a_1), &\ldots, \tilde{S}_n^{[0]}(a_r)) \\
&= \sum_{j=0}^{r} \sum_{I \subset \{1,\ldots,r\}, |I|=j} b_n^{-3j/2} \mathcal{K}_0(\{a_i' \tilde{S}_n : i \in I^c\}; \{a_i' \tilde{\xi}_{1n} : i \in I\}) \\
&= \sum_{j=0}^{s_r} \sum_{I \subset \{1,\ldots,r\}, |I|=j} b_n^{-3j/2} \mathcal{K}_0(\{a_i' \tilde{S}_n : i \in I^c\}; \{a_i' \tilde{\xi}_{1n} : i \in I\}) \\
&\quad + O(b_n^{-3(s_r+1)/2}).
\end{aligned}$$

By Lemma 4.1(ii) of [L], the $j=0$ term is $O(b_n^{-[r-2]/2})$ under assumptions (A.1)(i) and (A.2)–(A.5). Hence, consider $1 \leq j \leq r$ and fix an $I \in \{1, \ldots, r\}$ with $|I| = j$. Then $b_n^{-3j/2}|\mathcal{K}_0(\{a_i' \tilde{S}_n : i \in I^c\}; \{a_i' \tilde{\xi}_{1n} : i \in I\})|$ is bounded above by a sum of finitely many terms of the form

$$(6.11) \quad b_n^{-3j/2} |\mathcal{K}_0(\{a_i' \tilde{S}_n : i \in I^c\}; \{[_3\tilde{W}_{1n} + {_3\tilde{W}_{b_{0n}n}}]^{\alpha_i} [_1\tilde{S}_n]^{\beta_i} : i \in I\})|,$$

where $|\alpha_i| = 1 = |\beta_i|$ for all $i \in I$. To derive an upper bound on these terms, we make use of some known results on the cumulants of polynomial functions of a given set of random variables, for example, as given in Brillinger (1981). Consider the two-dimensional array

$$(6.12) \quad \begin{array}{ll}
(1,1) & (1,2) \\
(2,1) & (2,2) \\
\cdots & \cdots \\
(j,1) & (j,2) \\
(j+1,1) & \\
\cdots & \\
(r,1) &
\end{array}$$

and a partition $P_1 \cup \cdots \cup P_M$ of its entries. Following Brillinger (1981), we say that the sets $P_k$ and $P_l$ of the partition *hook* if there exist $(i_1, j_1) \in P_k$ and $(i_2, j_2) \in P_l$ with $i_1 = i_2$ and we say that $P_k$ and $P_l$ *communicate* if there exists a sequence of sets $P_{k_1} = P_k, P_{k_2}, \ldots, P_{k_N} = P_l$ from the partition such that $P_{k_i}$ and $P_{k_{i+1}}$ hook for all $i$. A partition is said to be *indecomposable* if all sets in the partition communicate. Then, for any collection of random variables $\{X_{p,q}\}$'s, indexed by the elements in the array (6.12), Theorem 2.3.2 of Brillinger (1981) yields

$$\mathcal{K}_0\left(\left\{\prod_{k=1}^{2} X_{i,k} : i = 1, \ldots, j\right\}; \{X_{i,1} : i = j+1, \ldots, r\}\right)$$

(6.13)
$$= \sum{}^{*} \mathcal{K}_0(\{X_{p,q} : (p,q) \in P_1\}) \cdots \mathcal{K}_0(\{X_{p,q} : (p,q) \in P_M\}),$$

where the sum $\sum^*$ above extends over all indecomposable partitions $P_1 \cup \cdots \cup P_M$ of the array (6.12). Note that if a partition of (6.12) is indecomposable, then all elements from the last $(r-j)$ rows must belong to a single set in the partition. W.l.g., we suppose that $P_1$ contains all the elements from these rows, i.e., $\{(i,1) : i = j+1, \ldots, r\} \subset P_1$. Thus, by (6.11) and (6.13), it is now enough to obtain a bound on terms of the form

(6.14)
$$b_n^{-3j/2} \mathcal{K}_0(\{a'_{1i}[_3\tilde{W}_{1n} + {}_3\tilde{W}_{b_0 n n}] : i \in P_{11}\};$$
$$\{a'_{2i}[_1\tilde{S}_n] : i \in P_{12}\}; \{a'_{3i}\tilde{S}_n : i \in P_{13}\})$$
$$\times \prod_{k=2}^{M} \mathcal{K}_0(\{x'_{ki}[_3\tilde{W}_{1n} + {}_3\tilde{W}_{b_0 n n}] : i \in P_{k1}\}; \{y'_{ki}[_1\tilde{S}_n] : i \in P_{k2}\})$$

for nonrandom vectors $a_{pi}$'s, $x_{ki}$'s and $y_{ki}$'s with vector-norm less than or equal to 1, and for indecomposable partitions $P_1 \cup \cdots \cup P_M$ of (6.12), where $P_1 = P_{11} \cup P_{12} \cup P_{13}$ with $|P_{13}| = r - j$ and $P_k = P_{k1} \cup P_{k2}$ for $k = 2, \ldots, M$. Further, $\sum_{k=2}^{M} |P_k| + \{|P_{11}| + |P_{12}|\} = 2j$, the total number of terms in the first $j$ rows and $\sum_{k=1}^{M} |P_{k1}| = j$, the total number of $[_3\tilde{W}_{1n} + {}_3\tilde{W}_{b_0 n n}]$ terms in (6.11), and $|P_k| \geq 2$ for all $k = 2, \ldots, M$. Under (A.2)–(A.5), by arguments in the proof of Lemma 4.1 of [L],

(6.15)
$$|\mathcal{K}_0(\{a'_{1i}[_3\tilde{W}_{1n} + {}_3\tilde{W}_{b_0 n n}] : i \in P_{11}\};$$
$$\{a'_{2i}[_1\tilde{S}_n] : i \in P_{12}\}; \{a'_{3i}\tilde{S}_n : i \in P_{13}\})|$$
$$= O(b_n^{-[|P_{12}| + |P_{13}| - 2]/2})$$

and for $k = 2, \ldots, M$,

(6.16)
$$|\mathcal{K}_0(\{x'_{ki}[_3\tilde{W}_{1n} + {}_3\tilde{W}_{b_0 n n}] : i \in P_{k1}\}; \{y'_{ki}[_1\tilde{S}_n] : i \in P_{k2}\})|$$
$$= O(b_n^{-(|P_{k2}| - 2)_+/2}).$$

Hence, it follows that the product in (6.14) is $O(b_n^{-3j/2} \times b_n^{-[|P_{12}| + |P_{13}| - 2]/2} \times \prod_{k=2}^{M} b_n^{-(|P_{k2}| - 2)_+/2})$. Using the conditions on $P_{kp}$'s, it can be shown that the largest order of these terms is attained when $|P_k| = 2$ for at most $\lfloor (2j-1)/2 \rfloor$-many $k$'s and $|P_1| \in \{r - j + 1, r - j + 2\}$ (so that $P_1$ contains

either one or two elements from the first $j$ rows). (The term corresponding to "$|P_1| = r - j + 1$ and $|P_k| = 2$ for all $k \geq 2$" actually does not appear in the sum $\sum^*$ due to the indecomposability condition.) Thus, the product in (6.14) is $O(b_n^{-3j/2} \times b_n^{-(r-j-2)/2}) = o(b_n^{-[r-2]/2})$ for all $j = 1, \ldots, s_r$. This completes the proof of the lemma. $\square$

For the next lemma, set $\tilde{V}_n(I) \equiv \prod_{j \in I} \prod_{p=1}^{r_j} \tilde{V}_{jn}(a_{jp})$ and $\tilde{W}_n(I) \equiv \prod_{j \in I} \prod_{p=1}^{r_j} a'_{jp} \tilde{W}_{jn}$, $I \subset \{1, \ldots, b_n\}$ where $a_{jp} \in \mathbb{R}^{d_0}$ with $\|a_{jp}\| \leq 1, r_j \in \mathbb{N}$. Also, for $r \geq 0$, $m \geq 3$, and $J \subset \{1, \ldots, b_{0n}\}$, let

$$H_{n,J}(t) = E \exp\left(\iota t' \sum_{j \in J} \tilde{W}_{jn}\right) \quad \text{and}$$

(6.17)
$$S_J^{(r)} \equiv S_J^{(r)}(t; m) = \iota b_n^{-1/2} \sum_{j \in J_r} \tilde{V}_{jn}(t),$$

$t \in \mathbb{R}^{d_0}$, where $J_r \equiv \{j : |j - i| > mr \text{ for all } i \in J\}$.

LEMMA 6.4. *Let $I \subset \{1, \ldots, b_{0n}\}$ with $1 \leq |I| \leq s + d_0 + 2$ and $\mathrm{DIAM}(I) \leq m_1$ for some $1 \leq m_1 \leq C[\log n]^2$. If assumptions (A.1)(i) and (A.2)–(A.5) hold, then*

(6.18) $\quad |E_t \tilde{V}_n(I)| \leq C[E|\tilde{W}_n(I)| \cdot |\theta_{n,I}(t)| + 2^{-K}] |H_n^{[0]}(t)|^{-1}$

*for all $\|t\| \leq e_{1n}$, where $e_{1n} = [\log n \log \log n]^{1/2}$, $K \equiv K_n = \lceil (\log n)^{3/2} \rceil$ and $\theta_{n,I}(t) = \sup\{|H_{n,J}(t + (x' \operatorname{SVEC}^{-1}(t_2), 0')')| : \|x\| \leq 2c_n/b_n^{3/2}, J \in \mathcal{J}_n\}$, with $\mathcal{J}_n$ being the collection of all subsets of $\{1, \ldots, b_{0n}\} \setminus I$ of size greater than $[b_{0n} - m_1 - (K+1)m]$, $m = 3\lceil K \log \log n \rceil + 1$.*

PROOF. Let $m_{0n} = \lceil K \log \log n \rceil$. Then, the following three cases are possible:

$$\text{(I) } I \cap \{1\}^{m_{0n}} \neq \varnothing; \quad \text{(II) } I \cap \{b_{0n}\}^{m_{0n}} \neq \varnothing;$$
$$\text{(III) } I \cap [\{1\}^{m_{0n}} \cup \{b_{0n}\}^{m_{0n}}] = \varnothing,$$

where $\{j\}^m = \{i : |i - j| \leq m\}$. We begin with case I. Note that in this case, $I \subset \{1, \ldots, m_{0n} + m_1 + 1\}$. With $J = \{1, \ldots, m_{0n} + m_1 + 1\} \cup \{b_n\}$ and $m = 3m_{0n} + 1$, define $\Delta_{1,r} = [\exp(S_J^{(r-1)} - S_J^{(r)}) - 1], r \geq 1$ and $\tilde{V}_{1n}(J) = \tilde{V}_n(I) \exp(\iota b_n^{-1/2} \sum_{j \in J} \tilde{V}_{jn}(t))$. Then using the iterative method of Tikhomirov (1980) as applied in the proof of Lemma 3.16 of [GH], one can show that

$$|H_n^{[0]}(t)| |E_t \tilde{V}_n(I)| = |E[\tilde{V}_{1n}(J) \exp(S_J^{(0)})]|$$

$$
(6.19) \quad \leq \sum_{r=1}^{K}\left|E\tilde{V}_{1n}(J)\left(\prod_{j=1}^{r-1}\Delta_{1,j}\right)\exp(S_J^{(r)})\right|
$$

$$
+\left|E\tilde{V}_{1n}(J)\left(\prod_{j=1}^{K}\Delta_{1,j}\right)\exp(S_J^{(K)})\right|.
$$

Note that for any $1 \leq j_0 < j_1 \leq b$, the variables $\{\tilde{W}_{jn}^{[0]} : j \leq j_0\}$ and $\{\tilde{W}_{jn}^{[0]} : j \geq j_1\}$ are functions $\{X_j : j \leq j_0\ell + \ell\}$ and $\{X_j : j \geq j_1\ell + 1\}$, respectively. Now, approximating $X_j$'s by $X_{j,m_{0n}}^{\dagger}$'s, we see that the corresponding variables $\{\tilde{W}_{jn}^{[0]\dagger} : j \leq j_0\}$ and $\{\tilde{W}_{jn}^{[0]\dagger} : j \geq j_1\}$ are measurable with respect to the $\sigma$-fields $\mathcal{D}_{-\infty}^{j_0\ell+\ell+m_{0n}}$ and $\mathcal{D}_{j_1\ell-m_{0n}+1}^{\infty}$, respectively. Next, approximate the variables $\tilde{V}_{1n}(J)$, $\Delta_{1,j}$'s and $S_J^{(r)}$'s using $X_{j,m_{0n}}$'s; call the approximating variables $\tilde{V}_{1n}^{\dagger}(J)$ $\Delta_{1,j}^{\dagger}$'s and $S_J^{(r)\dagger}$'s, respectively. Then, by (A.4), the $\alpha$-mixing co-efficient between $\sigma\langle\{\Delta_{1,j}^{\dagger}\}\rangle$ and $\sigma\langle\{\Delta_{1,j+2}^{\dagger}\}\rangle$ (and also, that between $\sigma\langle\{\Delta_{1,j}^{\dagger}\}\rangle$ and $\sigma\langle\{S_J^{(j+1)\dagger}\}\rangle$) is $O(\exp(-\kappa \times [(m-1)\ell - 2m_{0n}]))$. Hence, the last term in (6.19) is bounded above by

$$
(6.20) \quad \begin{aligned} & E|\tilde{V}(I)|\prod{}^{\prime}|\Delta_{1,j}|2^{K/2} \\ & \leq c_n^{r_0}2^{K/2}E\left(\prod{}^{\prime}|\Delta_{1,j}^{\dagger}|\right) + O(c_n^{r_0}K2^{K/2}\exp(-\kappa \cdot m_{0n})), \end{aligned}
$$

where $r_0 = \sum_{j\in I} r_j$ and where $\prod^{\prime}$ extends over all even indices $j \in \{1, \ldots, K\}$. Next, note that by (A.2) and Lemma 6.2, there exist independent random vectors $\tilde{W}_{1n}^{[0]*}$, $\tilde{W}_{b_{0n}n}^{[0]*}$, and $\{\tilde{W}_{jn}^{[0]*} : j \notin \{1\}^m \cup \{b_{0n}\}^m\} \equiv \{\tilde{W}_{jn}^{[0]\dagger} : j \notin \{1\}^m \cup \{b_{0n}\}^m\}$ such that with $u_n = \exp(-\kappa m_{0n}\ell)$,

$$
(6.21) \quad \begin{aligned} & P(\|\tilde{W}_{jn}^{[0]*} - \tilde{W}_{jn}^{[0]\dagger}\| > u_n) \\ & \leq C(d)(u_n^{-1}\exp(-\kappa[(m-1)\ell - 2m_{0n}])E\|\tilde{W}_{jn}^{[0]\dagger}\|)^{1/2} \\ & = O(\exp(-C(\kappa)m_{0n}\ell)) \end{aligned}
$$

for $j \in \{1, b_{0n}\}$. Next, define $\Delta_{1,r}^*$ by replacing $\tilde{W}_{jn}^{[0]\dagger}$'s in $\Delta_{1,r}^{\dagger}$ with $\tilde{W}_{jn}^{[0]*}$'s. Note that for any function $f : \mathbb{R}^{d_0} \times \mathbb{R}^k \to \mathbb{C}$ ($k \in \mathbb{N}$) with $|f(x;z) - f(y;z)| \leq C_1\|x-y\|$ and $|f(x;z)| \leq C_1$ for all $x \in \mathbb{R}^{d_0}, z \in \mathbb{R}^k$ and for some $C_1 \in (0,\infty)$,

$$
\sup_{z \in \mathbb{R}^k}|f(\tilde{W}_{jn}^{[0]\dagger};z) - f(\tilde{W}_{jn}^{[0]*};z)|
$$

$$
= \sup_{z \in \mathbb{R}^k}\{|f(\tilde{W}_{jn}^{[0]\dagger};z) - f(\tilde{W}_{jn}^{[0]*};z)|\mathbb{1}(\|\tilde{W}_{jn}^{[0]*} - \tilde{W}_{jn}^{[0]\dagger}\| \leq u_n)
$$

(6.22)
$$+ |f(\tilde{W}_{jn}^{[0]\dagger}; z) - f(\tilde{W}_{jn}^{[0]*}; z)| \mathbb{1}(\|\tilde{W}_{jn}^{[0]*} - \tilde{W}_{jn}^{[0]\dagger}\| > u_n)\}$$
$$\leq C_1 u_n + C_1 P(\|\tilde{W}_{jn}^{[0]*} - \tilde{W}_{jn}^{[0]\dagger}\| > u_n)$$
$$\leq C_1[C(d,\kappa)\exp(-C(\kappa)m_{0n}\ell)]$$

for $j \in \{1, b_{0n}\}$. Thus, by two applications of (6.22), we have $E|\Delta_{1,r}^\dagger - \Delta_{1,r}^*| = O(e_{1n}\exp(-C(\kappa)m_{0n}\ell))$ uniformly in $t \in B_{1n} \equiv \{x : \|x\| \leq e_{1n}\}$ and in $r \in \{1,\ldots,K\}$. Hence, using the independence of the $W_{1n}^{[0]*}$ and $W_{b_{0n}n}^{[0]*}$ with the rest of the $W_{jn}^{[0]*}$'s, the arguments in the proof of Lemma 3.16 of [GH] (for the second inequality below) and the above bound, uniformly in $t \in B_{1n}$, we have

|The last term in (6.19)|
$$\leq c_n^{r_0} 2^{K/2} E\left(\prod{}^{'} |\Delta_{1,j}^{[0]*}|\right)$$
$$+ O(c_n^{r_0} e_{1n} K 2^{K/2}[\exp(-C(\kappa)m_{0n}\ell) + \exp(-\kappa m_{0n})])$$

(6.23)
$$= c_n^{r_0} 2^{K/2} E\left\{E\left(\prod{}^{'} |\Delta_{1,j}^{[0]*}| \Big| \tilde{W}_{1n}^{[0]*}, \tilde{W}_{b_{0n}n}^{[0]*}\right)\right\} + O(\exp(-C(\kappa)m_{0n}))$$
$$\leq c_n^{r_0} 2^{K/2} E\left\{\prod{}^{'} E(|\Delta_{1,j}^{[0]*}| \big| \tilde{W}_{1n}^{[0]*}, \tilde{W}_{b_{0n}n}^{[0]*})\right\}$$
$$+ O(\exp(-C(\kappa)m_{0n}))$$
$$\leq c_n^{r_0} 2^{K/2} [C(\|\Sigma\|)\|t\| m^{1/2} b_n^{-1/2}]^{K/2} + O(\exp(-C(\kappa)m_{0n}))$$
$$\leq C(r_0, \|\Sigma\|) 2^{-K} + O(\exp(-C(\kappa)m_{0n})),$$

where in the last inequality, with $J_r = \{j : |j - i| > mr \text{ for all } i \in J\}$, we have used the bound

(6.24)
$$E(|\Delta_{1,r}^{[0]*}| \big| \tilde{W}_{1n}^{[0]*}, \tilde{W}_{b_{0n}n}^{[0]*})$$
$$\leq \|t\| b_n^{-1/2} \left[E\Big\|\sum_{j \in J_{r-1} \setminus J_r} \tilde{W}_{jn}^{[0]\dagger}\Big\|\right.$$
$$\left. + b_n^{-3/2} c_n E\Big\|\sum_{j \in J_{r-1} \setminus J_r} \tilde{W}_{jn}^{[0]\dagger}\Big\|\right]$$
$$\leq \|t\| b_n^{-1/2} C(\|\Sigma\|) m^{1/2}.$$

Next, consider the $r$th summand in the first term in (6.19) ($1 \leq K$). Note that
$$S_J^{(r)} = \iota t' \sum_{j \in J_r} \tilde{W}_{jn} + \iota b_n^{-3/2} [_3\tilde{W}_{1n} + {_3\tilde{W}_{b_{0n}n}}]' \text{SVEC}^{-1}(t_2) \sum_{j \in J_r} [_1\tilde{W}_{jn}],$$
where $t = (t_1', t_2')'$ with $t_2 \in \mathbb{R}^{\Lambda_0}$. Use Lemma 6.2 to construct independent copies of the random vectors $\{\tilde{W}_{jn}^{[0]\dagger} : j \leq m_{0n} + m_1 + 1 + (r-1)m\}$, $\{\tilde{W}_{jn}^{[0]\dagger} : b_{0n} - (r-1)m \leq j \leq b_{0n}\}$, and $\{\tilde{W}_{jn}^{[0]\dagger} : m_{0n} + m_1 + 1 + (rm+1) \leq j \leq b_{0n} - (rm+1)\}$. Then approximating $\tilde{W}_{jn}^{[0]}$'s with $\tilde{W}_{jn}^{[0]\dagger}$'s first and then $\tilde{W}_{jn}^{[0]\dagger}$'s with their independent copies $\tilde{W}_{jn}^{[0]*}$'s and using (6.21) and (6.22), uniformly in $t \in B_{1n}$, we get

$$\sum_{r=1}^{K} \left| E\tilde{V}_{1n}(J) \left(\prod_{j=1}^{r-1} \Delta_{1,j}\right) \exp(S_J^{(r)}) \right|$$

$$= \sum_{r=1}^{K} \left| E\tilde{V}_{1n}^*(J) \left(\prod_{j=1}^{r-1} \Delta_{1,j}^*\right) \exp(S_J^{(r)*}) \right|$$

$$+ O(c_n^{r_0} a_{0n} K^2 2^K [\exp(-C(\kappa) m_{0n} \ell) + \exp(-\kappa \cdot m_{0n})])$$

$$= \sum_{r=1}^{K} \left| E\tilde{V}_{1n}^*(J) \left(\prod_{j=1}^{r-1} \Delta_{1,j}^*\right) E\{\exp(S_J^{(r)*}) | W_{1n}^{[0]*}, W_{b_{0n}n}^{[0]*}\} \right|$$

(6.25) $\quad + O(\exp(-C(\kappa) m_{0n}))$

$$\leq \sum_{r=1}^{K} E \left| \tilde{V}_n^*(I) \left(\prod_{j=1}^{r-1} \Delta_{1,j}^*\right) \right|$$

$$\times \sup_{\|x\| \leq 2c_n/b_n} |H_{n,J_r}(t + (0', x' \text{SVEC}^{-1}(t_2))')|$$

$$+ O(\exp(-C(\kappa) m_{0n})).$$

Next, deleting the odd-$j$ terms in the product as in (6.20), using the weak dependence of the alternate sums as in (6.23), and the bound (6.24), one can show that uniformly over $t \in B_{1n}$,

$$\sum_{r=1}^{K} E \left| \tilde{V}_n^*(I) \left(\prod_{j=1}^{r-1} \Delta_{1,j}^*\right) \right|$$

$$\leq CE|\tilde{V}_n^*(I)| \sum_{r=1}^{K} 2^{r/2} [C(\|\Sigma\|)\|t\| m^{1/2} b_n^{-1/2}]^{r/2}$$

(6.26)
$$+ O(\exp(-C(\kappa)m_{0n}))$$
$$\leq CE|\tilde{V}_n^*(I)| + O(\exp(-C(\kappa)m_{0n})).$$

By (6.23), (6.25), (6.26), Lemma 6.4 is proved for case I.

The proofs of cases II and III are similar, with the set $J$ in the two cases now to be chosen as $J = \{1\} \cup \{b_{0n} - m_{0n} - m_1, \ldots, b_{0n}\}$ and $J = I \cup \{1\} \cup \{b_{0n}\}$, respectively. We omit the routine details. □

LEMMA 6.5. *Suppose that* (A.1)–(A.5) *hold. Let* $I_1, I_2, \subset \{1, \ldots, b\}$ *with* $\min\{I_2\} - \max\{I_1\} \geq m_1$ *for some integer* $C[\log n]^2 \leq m_1 \leq b_{0n} - |I_1| - |I_2|$. *Then there exists a constant* $C(\gamma) \in (0, \infty)$ *such that for all* $\|t\| \leq e_{1n}$, *and* $n \geq 1$,

(6.27) $\quad |E_t \tilde{V}_n(I_1)\tilde{V}_n(I_2) - E_t\tilde{V}_n(I_1)E_t\tilde{V}_n(I_2)| \leq C(\gamma) 2^{-K} |H_n^{[0]}(t)|^{-2},$

*where the variables* $\tilde{V}_n(I_k)$'*s are as defined above Lemma 5.4,* $\gamma = \sum_{k=1}^{2} \sum_{j \in I_k} r_j$, *and* $K = \lceil (\log n)^{3/2} \rceil$.

PROOF. Let $\{W_{jn}^{[0]**} : j = 1, \ldots, b_{0n}\}$ be an independent copy of $\{\tilde{W}_{jn}^{[0]} : j = 1, \ldots, b_{0n}\}$ and for a random vector $U = f(W_{jn}^{[0]} : j = 1, \ldots, b_{0n})$, let $U^{**} = f(W_{jn}^{[0]**} : j = 1, \ldots, b_{0n})$. Let $i_0 = \min\{I_2\}$. Define the variables $U_1 = \tilde{V}_n(I_1) - \tilde{V}_n^{**}(I_1)$, $U_2 = \tilde{V}_n(I_2) \prod_{j=i_0}^{b_{0n}} \exp(\iota b_n^{-1/2}[\tilde{V}_{jn}(t) + \tilde{V}_{jn}(t)^{**}])$ and $U_3 = \prod_{1 \leq j < i_0} \exp(\iota b_n^{-1/2}[\tilde{V}_{jn}(t) + \tilde{V}_{jn}^{**}(t)])$. Then it is easy to check that

$$|H_n^{[0]}(t)|^2 \times |E_t\tilde{V}_n(I_1)\tilde{V}_n(I_2) - E_t\tilde{V}_n(I_1)E_t\tilde{V}_n(I_2)| = |EU_1 U_2 U_3|.$$

Like in (6.17), write $S_n(r) = \iota b_n^{-1/2} \sum_{1 \leq j < i_0 - mr}[\tilde{V}_{jn} + \tilde{V}_{jn}^{**}]$, $\Delta_{2,r} = [\exp([S_n(r-1) - S_n(r)]) - 1]$, $0 \leq r \leq K$. Then, using Tikhomirov (1980)'s iterative method as in (6.19) (see also, Lemma 3.1 of [L] and Lemma 3.17 of [GH]), we get

(6.28)
$$EU_1U_2U_3 = \sum_{r=1}^{K} E\left\{U_1U_2\left(\prod_{j=1}^{r-1} \Delta_{2,j}\right) \exp(S_n(r))\right\}$$
$$+ E\left\{U_1U_2\left(\prod_{j=1}^{K} \Delta_{2,j}\right) \exp(S_n(K))\right\}.$$

Let $m_{0n} = \lceil \log n \log \log n \rceil$ and let $\mathbf{W}_{jn}^{[0]\dagger}$ be the $[2d_0] \times 1$ random vector with components $\tilde{W}_{jn,m_{0n}}^{[0]\dagger}$ and $\tilde{W}_{jn,m_{0n}}^{[0]**\dagger}$, $j = 1, \ldots, b_{0n}$. Next, consider the $r$th term in the sum in (6.28), $1 \leq r \leq K$. Note that by (A.4), the strong mixing coefficient between the variables $\{\mathbf{W}_{jn}^{[0]\dagger}: j \leq i_0 - mr\}$ and $\{\mathbf{W}_{jn}^{[0]\dagger}: j \geq i_0 - m(r-1)\}$ is bounded above by $2\kappa^{-1}\exp(-\kappa[(m-1)\ell - 2m_{0n}])$. Hence, using (A.2)(ii) and Lemma 5.2, we can construct independent copies $\{\mathbf{W}_{jn}^{[0]*}: j < i_0 - mr\}$ and $\{\mathbf{W}_{jn}^{[0]*}: j \geq i_0 - m(r-1)\}$ of these two sets of variables such that [cf. (6.21)]

$$P\Big(\max_{j \geq i_0 - m(r-1)} \|\mathbf{W}_{jn}^{[0]*} - \mathbf{W}_{jn}^{[0]\dagger}\| > u_n\Big)$$

(6.29)
$$\leq nC(d)(u_n^{-1}\exp(-\kappa[(m-1)\ell - 2m_{0n}])E\|\tilde{W}_{jn}^{[0]\dagger}\|)^{1/2}$$
$$= O(\exp(-C(\kappa)m_{0n}\ell)).$$

Hence, with the "natural" definition of the relevant $[\cdot]^\dagger$- and $[\cdot]^*$-variables (defined in terms of $\mathbf{W}_{jn}^{[0]\dagger}$'s and $\mathbf{W}_{jn}^{[0]*}$'s, respectively), uniformly in $1 \leq r \leq K$ and $t \in B_{1n}$, we have

$$E\Bigg\{U_1 U_2 \Bigg(\prod_{j=1}^{r-1} \Delta_{2,j}\Bigg) \exp(S_n(r))\Bigg\}$$
$$= E\Bigg\{U_1^\dagger U_2^\dagger \Bigg(\prod_{j=1}^{r-1} \Delta_{2,j}^\dagger\Bigg) \exp(S_n^\dagger(r))\Bigg\} + O(c_n^\gamma K 2^K \exp(-\kappa \cdot m_{0n}))$$
$$= E\Bigg\{U_1^* U_2^* \Bigg(\prod_{j=1}^{r-1} \Delta_{2,j}^*\Bigg) \exp(S_n^*(r))\Bigg\}$$
$$\quad + O(c_n^\gamma K 2^K [\exp(-\kappa \cdot m_{0n}) + \exp(-C(\kappa)m_{0n}\ell)])$$
$$= E\Bigg[E\Bigg(U_1^* U_2^* \Bigg(\prod_{j=1}^{r-1} \Delta_{2,j}^*\Bigg) \exp(S_n^*(r)) | \{\mathbf{W}_{jn}^{[0]*}: j > i_0 - [r-1]m\}\Bigg)\Bigg]$$
$$\quad + O(\exp(-C(\kappa)m_{0n}))$$
$$= O(\exp(-C(\kappa)m_{0n})),$$

since the conditional expectation in the first term is zero due to (i) the definition of $U_1^*$ and (ii) the symmetry of $U_1^* U_2^* [\prod_{j=1}^{r-1} \Delta_{2,j}^*] \exp(S_n^*(r))$ in the two $d_0 \times 1$ *i.i.d. components* of $\mathbf{W}_{jn}^*$ for all $j$. Next, using arguments similar to the proof of (6.23), one can show that the last term in (6.28) is $O(\exp(-C(\kappa)m_{0n}))$. Hence, (6.27) follows. □

LEMMA 6.6. *Suppose that assumptions* (A.1)–(A.5) *hold. Then, for any* $I \subset \{1, \ldots, b_{0n}\}$ *with* $|I| \leq C$, *and for any* $3 \leq m \leq b_{0n}/C$, *there exists* $L \geq C[b/m]$ *such that*

$$|H_n^{[0]}(t)||E_t\tilde{V}_n(I)| \leq Cc_n^\gamma[\beta_{1n}(t)]^{L-2} + C(|I|)Lc_n^\gamma[1 + \|t\|\ell m]\exp(-C(\kappa)m\ell)$$

for all $t \in \mathbb{R}^d$, where $\beta_{1n}(t) = \max\{E|E(\exp(\iota b_n^{-1/2} \sum_{|j-j_0| \le m} \tilde{V}_{jn}(t))|\tilde{\mathcal{D}}_{j_0})| : m < j_0 < b - m\}$, $\gamma = \sum_{j \in I} r_j$ and $\tilde{\mathcal{D}}_{j_0}$ is as in assumption (C.6). Further,

$$(6.30) \quad \beta_{1n}(t) \le \begin{cases} \exp(-\kappa m \|t\|^2/[8b_n]), & \text{for } \|t\| \le Cb_n/m^2, \\ \exp(-C(\kappa)(d_n b_n^{1/2})^{-2}\|t\|^2), & \text{for all } \|t\| \le 2\kappa b_n^{1/2} d_n. \end{cases}$$

PROOF. Let $I = \{j_1, \ldots, j_r\}$, $J_{0n} = \{1, \ldots, b_{0n}\}$, and $J_{1n} = \{j \in J_{0n} : |j - j_k| \ge 2m + 1 \text{ for all } 1 \le k \le r\}$. Next, define the indices $j_1^0, j_2^0, \ldots$ by

$$j_1^0 = \inf J_{1n},$$
$$j_{p+1}^0 = \inf\{j \in J_{1n} : j \ge j_p^0 + 7m\}, \qquad p = 1, 2, \ldots, L - 1,$$

where $L$ is the first integer $p$ for which the infimum is taken over an empty set. Then, it follows from the definitions of $j_1^0, \ldots, j_k^0$ that:

(i) $\lfloor (b_{0n} - r(2m))/7m \rfloor \le K \le \lceil b_{0n}/7m \rceil$,
(ii) $[j_p^0 - m, j_p^0 + m] \subset J_{0n}, p = 2, \ldots, K - 1$,
(iii) $j_{p+1}^0 \ge j_p^0 + 7m, p = 1, \ldots, K$, and
(iv) $(I)^m \cap \{j_1^0, \ldots, j_K^0\} = \varnothing$,

where $(I)^m = \{j \in J_{0n} : \text{there exists } i \in I \text{ with } |j - i| \le m\}$. Next, define the variables $A_p, B_p, p = 2, \ldots, L - 1$ and the variable $R$ by

$$A_p = \exp\left(\iota b_n^{-1/2} \sum_{|j - j_p^0| \le m} \tilde{V}_{jn}(t)\right),$$

$$B_p = \exp\left(\iota b_n^{-1/2} \sum_{j_p^0 + m + 1 \le j \le j_{p+1}^0 - m - 1} \tilde{V}_{jn}(t)\right),$$

$$R = \tilde{V}_n(I) \exp\left(\iota b_n^{-1/2} \sum_{j \in J_{2n}} \tilde{V}_{jn}(t)\right),$$

where $J_{2n} = \{j \in J_{0n} : j < j_1^0 + m \text{ or } j \ge j_L^0 - m\}$. Then it follows that

$$H_n^{[0]}(t) E_t \tilde{V}_n(I) = E\left[\left(\prod_{p=2}^{L-1} A_p B_p\right) R\right].$$

Next, define the variables $A_p^\dagger, B_p^\dagger$ and $R^\dagger$ by approximating $X_j$'s by $X_{j,q}^\dagger$ with $q = m\ell$ for $A_p^\dagger$ and with $q = m\ell/12$ for $B_p^\dagger$ and $R^\dagger$. Then, by (A.3), for all $p = 2, \ldots, L - 1$,

$$E|A_p - A_p^\dagger| \le C \cdot b_n^{-1/2}[1 + \|t\|] \cdot \ell m \exp(-\kappa m\ell)$$
$$E|B_p - B_p^\dagger| + E|R - R_p^\dagger| \le C(|I|) \cdot c_n^\gamma[1 + \|t\|] \cdot \ell m \exp(-\kappa m\ell/[12]).$$

Hence,

$$\left| ER \prod_{p=2}^{K-1} A_p B_p - ER^\dagger \prod_{p=2}^{K-1} A_p^\dagger B_p^\dagger \right|$$

(6.31)
$$\leq C(\kappa, |I|) K c_n^\gamma \|t\| \ell m \exp(-c(\kappa) m \ell).$$

Note that for any nonnegative function $\psi(X, Y)$ of random vectors $X \in \mathbb{R}^k$ and $Y \in \mathbb{R}^p$ on $(\Omega, \mathcal{F}, P)$ and a $\sigma$-field $\mathcal{C} \subset \mathcal{F}$ containing $\sigma\langle X \rangle$, a version of $E[\psi(X, Y)|\mathcal{C}]$ is given by

$$E[\psi(X, Y)|\mathcal{C}](\omega) = \psi_\mathcal{C}(X(\omega), \omega), \qquad \omega \in \Omega,$$

where $\psi_\mathcal{C}(x, \cdot) = \int \psi(x, y) \mu_\mathcal{C}(\cdot; dy)$, $x \in \mathbb{R}^k$, and $\mu_\mathcal{C}(\cdot; \cdot)$ is the regular conditional probability distribution of $Y$ on $\mathbb{R}^p$ given $\mathcal{C}$ [cf. Chapter 12, Athreya and Lahiri (2006)]. Hence, for any two $\sigma$-fields $\mathcal{C}$ and $\mathcal{D}$ containing $\sigma\langle X \rangle$ and for $\psi : \mathbb{R}^k \times \mathbb{R}^p \to [0, M]$, $M \in (0, \infty)$, we have

$$E|E[\psi(X, Y)|\mathcal{C}] - E[\psi(X, Y)|\mathcal{D}]|$$

$$= \int |\psi_\mathcal{C}(X(\omega), \omega) - \psi_\mathcal{D}(X(\omega), \omega)| P(d\omega)$$

$$= \int \left| \int \psi(X(\omega), y) \mu_\mathcal{C}(\omega, dy) - \int \psi(X(\omega), y) \mu_\mathcal{D}(\omega, dy) \right| P(d\omega)$$

$$\leq \int \left| \int_0^M [\mu_\mathcal{C}(\omega, \psi(X(\omega), \cdot)^{-1}([t, \infty))) \right.$$

(6.32)
$$\left. - \mu_\mathcal{D}(\omega, \psi(X(\omega), \cdot)^{-1}([t, \infty)))] dt \right| P(d\omega)$$

$$\leq \int_0^M \int |[\mu_\mathcal{C}(\omega, \psi(X(\omega), \cdot)^{-1}([t, \infty)))$$

$$- \mu_\mathcal{D}(\omega, \psi(X(\omega), \cdot)^{-1}([t, \infty)))]| P(d\omega)\, dt$$

$$\leq M \cdot E\left[ \sup_{B \in \mathcal{B}(\mathbb{R}^p)} |P(Y \in B|\mathcal{C}) - P(Y \in B|\mathcal{D})| \right].$$

Next, for each $p = 2, \ldots, L - 1$, define

$$\tilde{\mathcal{D}}_p = \sigma\langle \{\mathcal{D}_j : j \in \mathbb{Z}, j \notin [c_p, d_p]\} \rangle$$

and

$$\mathcal{D}_p^* = \sigma\langle \{\mathcal{D}_j : j \in [a_p - m\ell, c_p) \cup (d_p, b_p + m\ell]\} \rangle,$$

where $a_p = (j_p^0 - m)\ell + 1 - m\ell$, $b_p = (j_p^0 + m)\ell + m\ell$, $c_p = (j_p^0 - \lfloor \frac{m}{2} \rfloor)\ell + 1$, $d_p = (j_p^0 + \lfloor \frac{m}{2} \rfloor + 1)\ell$. Note that for each $p = 2, \ldots, L - 1$, $\{\tilde{W}_{jn}^{[0]\dagger} : j = 1, b_{0n}\}$ are measurable w.r.t. both the $\sigma$-fields $\tilde{\mathcal{D}}_p$ and $\mathcal{D}_p^*$. Hence, using (6.32) repeatedly (4 times) and assumption (A.5), one can show that

$$(6.33) \quad \max_{p=2,\ldots,L-1} |E(A_p^\dagger|\tilde{\mathcal{D}}_p) - E(A_p^\dagger|\mathcal{D}_p^*)| \leq C\exp(-C(\kappa)m\ell).$$

Also, note that the variables $R^\dagger(\prod_{p=2}^{q-1} A_p^\dagger B_p^\dagger)$, $B_q^\dagger$, $\prod_{p=q+1}^{L-1} B_p^\dagger E(A_p^\dagger|\mathcal{D}_p^*)$ are all measurable w.r.t. $\tilde{\mathcal{D}}_q$ for every $2 \leq q \leq L-1$. Hence, it follows that

$$(6.34) \quad \begin{aligned} &\left| ER^\dagger \prod_{p=2}^{L-1} A_p^\dagger B_p^\dagger - ER^\dagger \prod_{p=2}^{L-1} B_p^\dagger E(A_p^\dagger|\mathcal{D}_p^*) \right| \\ &\leq \sum_{q=2}^{L-1} \left| ER^\dagger \left( \prod_{p=2}^{q-1} A_p^\dagger B_p^\dagger \right) B_q^\dagger [A_q^\dagger - E(A_q^\dagger|\tilde{\mathcal{D}}_q)] \prod_{p=q+1}^{L-1} B_p^\dagger E(A_p^\dagger|\mathcal{D}_p^*) \right| \\ &\quad + Cc_n^\gamma \sum_{q=2}^{L-1} E|E(A_q|\tilde{\mathcal{D}}_q) - E(A_q^\dagger|\mathcal{D}_q^*)| \\ &\leq C(\kappa)c_n^\gamma L\exp(-\kappa m\ell), \end{aligned}$$

since the first term vanishes. Now, using the fact that $\mathcal{D}_p^*$ and $\mathcal{D}_{p+1}^*$ are separated by a distance $\geq Cm\ell$ for all $p$, and using (6.31), (6.33), (6.34), one can conclude that

$$\begin{aligned} &|H_n^{[0]}(t) E_t \tilde{V}_n(I)| \\ &\leq Cc_n^\gamma E \left| \prod_{p=2}^{L-1} E(A_p^\dagger|\mathcal{D}_p^*) \right| + C(|I|)Lc_n^\gamma(1+\|t\|\ell m)\exp(-C(\kappa)m\ell) \\ &\leq Cc_n^\gamma \prod_{p=2}^{L-1} E|E(A_p|\tilde{\mathcal{D}}_p)| + C(|I|)Lc_n^\gamma(1+\|t\|\ell m)\exp(-C(\kappa)m\ell). \end{aligned}$$

The first part of Lemma 6.6 follows from this.

Next, we consider (6.30) and prove the first inequality. Let $Z$ be a zero mean $\mathbb{R}^k$-valued random vector on $(\Omega, \mathcal{F}, P)$ with $E\|Z\|^3 < \infty$ and let $\mathcal{C} \subset \mathcal{F}$ be a $\sigma$-field. Denote the conditional distribution of $Z$ given $\mathcal{C}$ on $(\mathbb{R}^k, \mathcal{B}(\mathbb{R}^k))$ by $\mu_{\mathcal{C}}(\cdot;\cdot)$. Then

$$\begin{aligned} &|E\exp(\sqrt{-1}t'Z|\mathcal{C})|^2 \\ &= \iint \exp(\sqrt{-1}t'(x-y))\mu_{\mathcal{C}}(\cdot,dx)\mu_{\mathcal{C}}(\cdot,dy) \\ &\leq 1 - E((t'Z)^2|\mathcal{C}) + 2E(|t'Z|^3|\mathcal{C}) \qquad \text{a.s. } (P) \end{aligned}$$

for all $t \in \mathbb{R}^k$. Hence, by (A.2), it follows that for all $n > C$,

$$\left[ E \left| E \left\{ \exp\left(\iota \left[ b_n^{-1/2} \sum_{|j-j_0| \leq m} \tilde{V}_{jn}(t) \right] \right) \middle| \tilde{\mathcal{D}}_{j_0} \right\} \right| \right]^2$$

$$\leq 1 - b_1^{-1} E \left( \sum_{|j-j_0| \leq m} \tilde{V}_{jn}(t) \right)^2$$

$$+ 2 b_n^{-3/2} E \left| \sum_{|j-j_0| \leq m} \tilde{V}_{jn}(t) \right|^3$$

$$\leq \exp(-b_n^{-1}[\kappa/2]m\|t\|^2 + C b_n^{-3/2}\|t\|^3 m^3)$$

$$\leq \exp\left(-\frac{\kappa}{4} b_n^{-1} m \|t\|^2\right)$$

for all $\|t\| \leq C b_n^{1/2}/m^2$, which proves the first bound on $\beta_{1n}(\cdot)$ in (6.30).

Next, consider the second bound on $\beta_{1n}(\cdot)$ in (6.30). Note that by iterating the inequality "$1 - \cos(2x) \leq 4(1 - \cos x)$ for all $x \in \mathbb{R}$" for $r$-times ($r \in \mathbb{N}$), we have $1 - \cos(2^r x) \leq 4^r(1 - \cos x), x \in \mathbb{R}$. Hence, for any random variable $Z$ and any $\sigma$-field $\mathcal{C} \subset \mathcal{F}$, if the conditional distribution of $Z$ given $\mathcal{C}$ is symmetric about the origin, then

$$1 - E[\exp(i 2^r u Z)|\mathcal{C}] \leq 4^r(1 - E[\exp(iuZ)|\mathcal{C}])$$

for all $u \in \mathbb{R}$. Hence, for a random variable $Z_1$ and a subsigma-field $\mathcal{G}$ of $\mathcal{F}$, writing $\mu_{\mathcal{G}}(\cdot;\cdot)$ for the conditional distribution of $Z_1$ given $\mathcal{G}$, and employing the above inequality, one can conclude (cf. page 223, [GH]) that

$$1 - |E(\exp(i 2^r u Z_1)|\mathcal{G})|^2$$

$$= 1 - \iint \exp(i 2^r u [z_1 - z_2]) \mu_{\mathcal{G}}(\cdot, dz_1) \mu_{\mathcal{G}}(\cdot, dz_2)$$

$$\leq 4^r \left[ 1 - \iint \exp(iu[z_1 - z_2]) \mu_{\mathcal{G}}(\cdot, dz_1) \mu_{\mathcal{G}}(\cdot, dz_2) \right]$$

$$= 4^r [1 - |E(\exp(iuZ_1)|\mathcal{G})|^2]$$

for all $u \in \mathbb{R}$. Hence, it follows that

$$1 - \left| E \left( \exp\left( \iota \sum_{|j-j_0| \leq m} \tilde{V}_{jn}(t) \right) \middle| \tilde{\mathcal{D}}_{j_0} \right) \right|^2$$

(6.35) $$\geq 4^{-r} \left( 1 - \left| E \left( \exp\left( \iota 2^r \sum_{|j-j_0| \leq m} \tilde{V}_{jn}(t) \right) \middle| \tilde{\mathcal{D}}_{j_0} \right) \right|^2 \right)$$

$$\geq 4^{-r}\left[1 - \sup_{\|t\|/\sqrt{2}\leq\|y\|\leq\sqrt{2}\|t\|}\left|E\left(\exp\left(\iota 2^r \sum_{|j-j_0|\leq m} y'\tilde{W}_{jn}\right)\Big|\tilde{\mathcal{D}}_{j_0}\right)\right|^2\right]$$

for all $t \in \mathbb{R}^{d_0}$ and $n > C$, where in the last step, we have used (6.1) to conclude that there exists a $C > 0$ such that for all $n > C$,

$$2^{-1}\|t\|^2 \leq \|(t'_1 + x' \operatorname{SVEC}^{-1}(t_2), t'_2)'\|^2$$
$$\leq 2\|t\|^2,$$

uniformly in $\|x\| \leq 2c_n b_n^{-3/2}$ and in $t = (t'_1, t'_2)' \in \mathbb{R}^d \times \mathbb{R}^{d_1}$.

Next, set $\kappa_1 = 2\kappa$ where $\kappa$ as in (A.6). Then, given a $t \in \mathbb{R}^{d_0}$ with $\|t\| \leq \kappa_1 d_n$, let $r \equiv r_t$ be such that

$$\frac{2^{r+1}\|t\|}{\sqrt{2}d_n} > \kappa_1 \geq \frac{2^r\|t\|}{\sqrt{2}d_n}.$$

Then, $4^{-r} \geq [\|t\|/\{\sqrt{2}\kappa_1 d_n\}]^2$ and by (A.6),

$$E\left[\sup_{\|t\|/\sqrt{2}\leq\|y\|\leq\sqrt{2}\|t\|}\left|E\left(\exp\left(\iota 2^r \sum_{|j-j_0|\leq m} y'\tilde{W}_{jn}\right)\Big|\tilde{\mathcal{D}}_{j_0}\right)\right|^2\right]$$

$$\leq E\left[\sup_{\kappa_1 d_n/2\leq\|x\|\leq 2\kappa_1 d_n}\left|E\left(\exp\left(\iota 2^r \sum_{|j-j_0|\leq m} x'\tilde{W}_{jn}\right)\Big|\tilde{\mathcal{D}}_{j_0}\right)\right|^2\right]$$

$$\leq E\left[\sup_{\|x\|\geq\kappa d_n}\left|E\left(\exp\left(\iota 2^r \sum_{|j-j_0|\leq m} x'W_{jn}\right)\Big|\tilde{\mathcal{D}}_{j_0}\right)\right|^2\right]$$

$$+ 2\sum_{j=j_0-m_n}^{j_0+m_n} P(\tilde{W}_{jn} \neq W_{jn})$$

$$\leq (1-\kappa) + 2\sum_{j=j_0-m_n}^{j_0+m_n} E\|W_{jn}\|^3 c_n^{-3}$$

$$\leq (1-[\kappa/2])$$

for all $n \geq C$. Hence, from (6.35), for all $\|t\| \leq 2\kappa d_n$,

$$E\left|E\left(\exp\left(\iota \sum_{|j-j_0|\leq m} \tilde{V}_{jn}(t)\right)\Big|\tilde{\mathcal{D}}_{j_0}\right)\right|$$

$$\leq \left\{E\left|E\left(\exp\left(\iota \sum_{|j-j_0|\leq m} \tilde{V}_{jn}(t)\right)\Big|\tilde{\mathcal{D}}_{j_0}\right)\right|^2\right\}^{1/2}$$

$$\leq \left(1 - \frac{\|t\|^2}{2\kappa_1^2 d_n^2} \cdot \frac{\kappa}{2}\right)^{1/2}$$

$$\leq \exp(-C(\kappa)\|t\|^2/d_n^2).$$

It is easy to see that this bound holds uniformly over $j_0 \in \{m+1, \ldots, b_{0n} - m\}$ and therefore, for $n > C$, we have

(6.36) $\quad \beta_{1n}(t) \leq \exp(-C(\kappa)(d_n b_n^{1/2})^{-2}\|t\|^2) \quad$ for all $\|t\| < 2\kappa b_n^{1/2} d_n$.

This completes the proof of (6.30) and hence of Lemma 6.6. $\square$

PROOF OF THEOREM 5.1. In view of the moment condition (A.2), by Lemma 4.1(i) of [L] and (the proof of) Lemma 6.3 above, and by Theorem 10.1, Corollary 11.5 and Lemma 11.6 of [BR] and by arguments similar to the proof of relation (20.39) in [BR], it is enough to show that

(6.37) $\quad \sup_{|\alpha| \leq d_0 + 1} \int_{\|\delta_n t\| \leq 1} |D^\alpha(H_n^{[0]}(t) - \tilde{\Upsilon}_{s,n}(t))| \, dt = O(\delta_n),$

where $\tilde{\Upsilon}_{s,n}(\cdot)$ is defined by replacing the cumulants of $S_n + b_n^{-3/2}\xi_{1n}$ in $\hat{\Upsilon}_{s,n}(\cdot)$ with those of the truncated and centered version $\tilde{S}_n + b_n^{-3/2}\tilde{\xi}_{1n}$.

First, consider the integral in (6.37) over the set $B_{1n} \equiv \{t \in \mathbb{R}^{d_0} : \|t\| \leq e_{1n}\}$, where $e_{1n} = [\log n][\log \log n]$. By Lemma 9.7 in [BR] and the arguments on pages 231 and 232 in the proof of Lemma 3.33 of [GH], it is enough to show that there exists an $\eta \in (0, 1)$ such that

$$\left|\frac{\partial^r}{\partial \varepsilon^r}\bigg|_{\varepsilon=0} R_n(t + \varepsilon a)\right| \leq C(\eta)(1 + \|t\|^{s+1})b_n^{-\eta - (s-2)/2}$$

for all $r = 0, \ldots, d_0 + 1$, $t \in A_{1n}$, where $R_n(t) = [\int_0^1 (1-u)^s \mathcal{K}_{ut}^{[0]}([\tilde{S}_n^{[0]}(\iota \times t)]^{\diamond(s+1)}) \, du]/s!$ is as in (6.7). In view of (6.3), (6.4), (6.7) and the multilinearity of the semi-invariants, it is enough to estimate the sums

(6.38) $\quad \sum_{i=0}^{b_{0n} - r} \sum^{(i,r)} |\mathcal{K}_t^{[0]}(V_{j_1 n}(a_1), \ldots, V_{j_r n}(a_r))|, \quad r = s+1, \ldots, s+d_0+1,$

for $t \in B_{1n}$, where $a_1, \ldots, a_r \in \mathbb{R}^{d_0}$ with $\|a_i\| \leq 1$ for all $1 \leq i \leq r$ and where the sum $\sum^{(i,r)}$ extends over all indices $1 \leq j_1 \leq \cdots \leq j_r$ with maximal gap $i$. Note that for any $s + 1 \leq r \leq s + d_0 + 1$,

(6.39)
$$E\|\tilde{W}_{jn}\|^r \leq c_n^{r-s-1}[E\|\tilde{W}_{jn}\|^{s+1}\mathbb{1}(\|\tilde{W}_{jn}\| \leq b_n^{1/4})$$
$$+ c_n E\|\tilde{W}_{jn}\|^s \mathbb{1}(\|\tilde{W}_{jn}\| > b_n^{1/4})]$$
$$\leq C(r,s)c_n^{r-s-[\kappa/4]},$$

uniformly in $1 \leq j \leq b_{0n}$.

Also, by relation (6.1), Lemma 6.4 and (A.2), for any $t = (t_1', t_2')' \in \mathbb{R}^d \times \mathbb{R}^{d_1}$, there exists a $u \in [-1, 1]$ such that

$$
\begin{aligned}
|H_n^{[0]}&(t) - H_n(t)| \\
&\leq |b_n^{-3/2} E(\iota t' \tilde{\xi}_{1n} \exp(\iota t' \tilde{S}_n + \iota u b_n^{-3/2} t' \tilde{\xi}_{1n}))| \\
&\leq b_n^{-2} \sum_{j=1}^{b_{0n}} |E\{([_3\tilde{W}_{1n}] + [_3\tilde{W}_{b_{0n}n}])' \\
&\qquad \times \mathrm{SVEC}^{-1}(t_2)[_1\tilde{W}_{jn}] \exp(\iota t' \tilde{S}_n + \iota u b_n^{-3/2} t' \tilde{\xi}_{1n})\}| \\
&\leq C b_n^{-2} \|t\| \sum_{j=1}^{b_{0n}} \max_{I=\{1,j\},\{j,b_{0n}\}} |E\tilde{W}_n(I) \exp(\iota t' \tilde{S}_n + \iota u b_n^{-3/2} t' \tilde{\xi}_{1n})| \\
&\leq C b_n^{-1} \|t\| \{\beta_{2n}(t) + 2^{-K}\},
\end{aligned}
\tag{6.40}
$$

where $\beta_{2n}(t) = \max\{|\theta_{n,I}(t)| : I = \{k, j\}, k \in \{1, b_{0n}\}, j \in \{1, \ldots, b_{0n}\}\}$ and $\theta_{n,I}(t)$ is as in Lemma 6.4. By arguments in the proofs of Lemma 3.33 of [GH83] and Lemma 4.4 of [L] (with $s = 3$, $q = 0$), we have

$$
\begin{aligned}
\sup\{\beta_{2n}(t)/|H_n(t)| : t \in B_{1n}\} &\leq C \quad \text{for } n > C \quad \text{and} \\
H_n(t) &= \exp(-t'\Sigma_n t/2 + O(\|b_n^{-1/2} t\|^3))
\end{aligned}
\tag{6.41}
$$

uniformly in $t \in B_{1n}$,

where $\Sigma_n = \mathrm{Var}(\tilde{S}_n)$.

Next, let $e_{2n} = (\log n)^2$ and fix $r \in \{s+1, \ldots, s+d_0+1\}$. Then, applying (6.39), (6.40), (6.41), and Lemma 6.4 for $i \leq e_{2n}$ and applying (6.6), (6.40), (6.41) and Lemma 6.5 for $e_{2n} \leq i \leq b_{0n} - 1$, from (6.38), we get

$$
\begin{aligned}
b_n^{-r/2} &\sum_{i=0}^{b_{0n}-r} \sum^{(i,r)} |\mathcal{K}_t^{[0]}(V_{j_1 n}(a_1), \ldots, V_{j_r n}(a_r))| \\
&\leq b_n^{-r/2+1} \sum_{i=0}^{e_{2n}} i^{r-1} \left[ \max_{|I|=r} |E\tilde{W}_n(I)| \cdot C + 2^{-K} / \exp(-C\|t\|^2) \right] \\
&\quad + b_n^{-r/2+1} \sum_{i=e_{2n}+1}^{b_{0n}-1} i^{r-1} [2^{-K} / \exp(-C\|t\|^2)] \\
&\leq C b_n^{-(r-2)/2} c_n^{r-s-[\kappa/4]} e_{2n}^r + C 2^{-K} \exp(C e_{1n}^2) \\
&\leq C b_n^{-(r-2)/2} c_n^{r-s-[\kappa/4]} (\log n)^C,
\end{aligned}
$$

uniformly in $t \in B_{1n}$, for $n > C$. Hence,

$$\max_{|\alpha| \leq d_0+1} \int_{t \in B_{1n}} |D^\alpha(H_n^{[0]}(t) - \hat{\Upsilon}_{s,n}(t))| \, dt = O(b_n^{-(s-2)/2}[\log n]^{-2}).$$

Also, by (A.2),

$$\max_{|\alpha| \leq d_0+1} \int_{\|t\| \geq e_{1n}} |D^\alpha \hat{\Upsilon}_{s,n}(t)| \, dt = O(b_n^{-a})$$

for any $a > 0$. Therefore, it remains to show that

(6.42) $$\max_{|\alpha| \leq d_0+1} \int_{\|t\| \geq e_{1n}} |D^\alpha H_n^{[0]}(t)| \, dt = O(b_n^{-(s-2)/2}[\log n]^{-2}).$$

Next, set $m = \log n$ if $\ell \leq n^{1/4}$ and set $m = C = C(s, d, \kappa)$ (a large but finite constant depending on $s, d, \kappa$), otherwise. Then, using the first bound of (6.30) for $e_{1n} \leq \|t\| \leq Cb_n^{1/2}/m^2$, the second bound of (6.30) for $b_n^{1/2}/m^2 \leq \|t\| \leq 2\kappa b_n^{1/2} d_n$, and Lemma 6.6 and condition (A.6) for $2\kappa b_n^{1/2} d_n \leq \|t\| \leq b_n^{(s-2)/2}[\log n]^2$, one gets (6.42). This completes the proof of the theorem. $\square$

### 6.3. Proof of the main result from Section 2.

PROOF OF PROPOSITION 2.1. As indicated in Section 2, we set $\mathcal{D}_j = \sigma\langle \varepsilon_j \rangle$, $j \in \mathbb{Z}$. Let $\mathcal{D}_{j_0}^* = \sigma\langle \varepsilon_j : j \neq j_0 \ell \rangle$, $j_0 \in J_n$, where $J_n$ is as in (A.6). Then, $\tilde{\mathcal{D}}_{j_0} \subset \mathcal{D}_{j_0}^*$, and hence, by the properties of the conditional expectation, it is enough to show that (2.3) holds with $\tilde{\mathcal{D}}_{j_0}$ replaced with $\mathcal{D}_{j_0}^*$. Using the independence of the $\varepsilon_i$'s, for any $t = (t_1', t_2')' \in \mathbb{R}^d \times \mathbb{R}^{d_1}$ and $j_0 \in J_n$, after some lengthy algebra, we get

(6.43)
$$\left| E\left( \exp\left( \iota t' \sum_{j=j_0-m}^{j_0+m} W_{jn} \right) \Big| \mathcal{D}_{j_0}^* \right) \right|$$
$$= \left| E\left\{ \exp\left( \iota \sum_{j=(j_0-m-1)\ell+1}^{(j_0+m)\ell} [\ell^{-1/2} t_1' X_j + \ell^{-1} t_2' Y_{jn}] \right) \Big| \{\varepsilon_j, j \neq j_0 \ell\} \right\} \right|$$
$$= |E\{\exp(\iota[\ell^{-1/2} t_1' L_n + \ell^{-1} t_2' R_n(j_0)]\varepsilon_{j_0\ell} + \varepsilon_{j_0\ell}' M_n(t_2)\varepsilon_{j_0\ell})\}|,$$

say,

where

$$L_n = \sum_{j=-(m+1)\ell}^{m\ell} A_j,$$

$$M_n(t_2) = \ell^{-1} \sum_{j=-(m+1)\ell}^{m\ell} \sum_{k=0}^{\ell} w_{kn}[A'_j \text{VEC}^{-1}(t_2)A_{j+k} + A'_{j+k}\text{VEC}^{-1}(t_2)A_j]$$

and $R_n(j_0)$ is a $\mathcal{D}^*_{j_0}$-measurable $d_1 \times d$ matrix-valued random element satisfying $E\|R_n(j_0)\|^2 \leq C\ell$. Here and in the following, the generic constants and the order symbols do not depend on $j_0 \in J_n$ (i.e., uniform over $j_0 \in J_n$).

Next, write $\tilde{\varepsilon} = \varepsilon_{j_0\ell}$ (for notational simplicity). By interchanging the order of summation in the second summand in $M_n(t_2)$ and using (6.1), one can show that

(6.44)
$$\varepsilon'_{j_0\ell} M_n(t_2) \varepsilon_{j_0\ell} = \tilde{\varepsilon}' \left[ \sum_{j=-m\ell+1}^{m\ell} A'_j \text{VEC}^{-1}(t_2) \left( \sum_{|k|\leq \ell} \check{w}_{kn} A_{j+k} \right) \right.$$
$$\left. + \sum^{(*)} A_j \text{VEC}^{-1}(t_2) \check{w}_{jkn} A_{j+k} \right] \tilde{\varepsilon}$$
$$= t'_2 \text{VEC} \left( \sum_{j=-m\ell+1}^{m\ell} \sum_{|k|\leq \ell} \check{w}_{kn}[A'_j \tilde{\varepsilon}][A_{j+k}\tilde{\varepsilon}]' \right.$$
$$\left. + \sum^{(*)} \check{w}_{jkn}[A'_j\tilde{\varepsilon}][A_{j+k}\tilde{\varepsilon}]' \right),$$

where $\check{w}_{kn} = 1$ if $k = 0$ and $\check{w}_{kn} = w_{kn}$ for all $1 \leq k \leq \ell$, and where the sum $\sum^{(*)}$ extends over all $j, k$ satisfying $j \in [-(m+1)\ell+1, -m\ell] \cup [m\ell+1, (m+1)\ell]$ and $1 \leq k \leq \ell$, and where $\check{w}_{jkn}$'s some real numbers satisfying $|\check{w}_{jkn}| \leq C$ for all $(j,k)$ under $\sum^{(*)}$. Note that the VEC($\cdot$)-term in (6.44) is a *linear combination* of the random vector $U \equiv \text{VEC}(\tilde{\varepsilon}\tilde{\varepsilon}')$; We write this terms as $D_n U$, where $D_n$ is a $d_1 \times d_1$ matrix. Next, using the condition $\lim_{n\to\infty} w_{kn} = 1$ for every $k \geq 1$ and the summability of $\|A_i\|$'s, one can show that for any $x \in \mathbb{R}^d$,

$$\left\| \left( \sum_{j=-m\ell+1}^{m\ell} \sum_{|k|\leq \ell} \tilde{w}_{kn}[A'_j x][A_{j+k}x]' \right.\right.$$
$$\left.\left. + \sum^{(*)} w_{jkn}[A'_j x][A_{j+k}x]' \right) - A_\infty xx' A_\infty \right\| \to 0$$

as $n \to \infty$. This implies that the matrix $D_n$ has a limit (say $D$), defined by

$$DU = \text{VEC}(A_\infty \tilde{\varepsilon}\tilde{\varepsilon}' A_\infty).$$

Next, let $Q(\cdot)$ and $Q_n(\cdot)$ denote the probability distributions of $V \equiv ([A_\infty \tilde{\varepsilon}]'; [DU]')'$ and $V_n \equiv ([L_n\tilde{\varepsilon}]', [D_n U]')'$, respectively. Since the distribution of $\tilde{\varepsilon}$ has a nontrivial absolutely continuous component, by Lemma 2.2 of Bhattacharya and Ghosh (1978), the $k$-fold convolution $Q^{*k}$ of $Q$ (where $k = d_0$) has an absolutely continuous component with respect to the Lebesgue measure on $\mathbb{R}^{d_0 k}$, with density $f^{(k)}$ (say). By Lusin's theorem [cf. Theorem 2.5.12, Athreya and Lahiri (2006)], without loss of generality, we may assume that $f^{(k)}$ is continuous. (In fact, the continuity of $f^{(k)}$ holds over a smaller set of *positive measure*, but that is enough for our purpose; we may just consider the absolutely continuous component restricted to this set.) Since $D_n \to D$, $L_n \to A_\infty$, and $Q_n^{*k}$ is obtained from $Q^{*k}$ by a linear transformation that is nonsingular for $n$ large, it follows that for large $n$, $Q_n^{*k}$ has an absolutely continuous component with respect to the Lebesgue measure on $\mathbb{R}^{d_0 k}$, with a density $f_n^{(k)}$ (say). Now using the standard transformation technique formula for the density $f_n^{(k)}$ in terms of $f^{(k)}$, we conclude that

$$f_n^{(k)} \to f^{(k)} \qquad \text{as } n \to \infty \text{ (a.e.)}.$$

Let $\varphi_n(t)$ denote the characteristic function of $V_n$ and write $Q_{n,s}^{*k}$ for the singular part of $Q_n^{*k}$. Then, it follows that for $t \in \mathbb{R}^{d_0}$,

$$|\varphi_n(t)|^k = \left| \int \exp(\iota t'(x_1 + \cdots + x_k)) \, dQ_n^{*k}(x) \right| \qquad [\text{where } x = (x_1', \ldots, x_k')']$$

$$\leq Q_{n,s}^{*k}(\mathbb{R}^{d_0 k}) + \left| \int \exp(\iota t'(x_1 + \cdots + x_k)) f_n^{(k)}(x) \, dx \right|$$

$$\leq Q_{n,s}^{*k}(\mathbb{R}^{d_0 k}) + \int |f_n^{(k)} - f^{(k)}| \, dx$$

$$+ \left| \int \exp(\iota t'(x_1 + \cdots + x_k)) f^{(k)}(x) \, dx \right|.$$

Hence, for every $\kappa \in (0, \infty)$, there exist a $\delta \in (0, 1/2)$ and $n_0 = n_0(\kappa)$ such that for all $n \geq n_0$,

(6.45) $$\sup_{\|t\| \geq \kappa} |\varphi_n(t)| \leq 1 - 2\delta.$$

Now using (6.43)–(6.45) and writing $B_{1n} = \{\|R_n(j_0)\| \leq C_1(\delta)\ell^{1/2}\}$, we have, for $n \geq n_0$,

$$E\left( \sup_{\|t\| \geq C\ell} \left| \left( \exp\left( \iota t' \sum_{j=j_0-m}^{j_0+m} W_{jn} \right) \Big| \mathcal{D}_{j_0}^* \right) \right| \right)$$

$$= E\left( \sup_{\|t\| \geq C\ell} |\varphi(\ell^{-1/2} t_1 + \ell^{-1} t_2' R_n(j_0), \ell^{-1} t_2)| \right) \qquad \text{where } t = (t_1', t_2')'$$

$$\leq \max\Big\{ E\Big( \sup_{\|t_2\|\geq \kappa\ell} |\varphi(\ell^{-1/2}t_1 + \ell^{-1}R_n(j_0)'t_2, \ell^{-1}t_2)|\Big),$$

$$E\Big( \sup_{\|t\|\geq C\ell, \|t_2\|\leq \kappa\ell} \{|\varphi(\ell^{-1/2}t_1 + \ell^{-1}R_n(j_0)'t_2,$$

$$\ell^{-1}t_2)|\mathbb{1}(B_{1n})\}\Big) + P(B_{1n}^c)\Big\}$$

$$\leq \max\Big\{ \sup_{\|s_2\|\geq \kappa, s_1\in\mathbb{R}^d} |\varphi(s_1, s_2)|,$$

$$\sup_{\|s_1\|\geq C(\kappa,\delta)\ell^{1/2}, s_2\in\mathbb{R}^{d_1}} |\varphi(s_1, s_2)| + \frac{E\|R_n(j_0)\|^2}{C_1(\delta)^2\ell} \Big\}$$

$$\leq 1 - \delta,$$

where in the step before the last one, we have used the fact that $\|t_2\| \leq \kappa\ell, \|t\| \geq C\ell \Rightarrow \|t_1\| \geq (C-\kappa)\ell$ which, in turn, implies that on the set $B_{1n}$, $\|\ell^{-1/2}t_1 - \ell^{-1}R_n(j_0)'t_2\| \geq \|\ell^{-1/2}t_1\| - \ell^{-1}\|R_n(j_0)\|\|t_2\| \geq (C-\kappa-C_1(\delta))\ell^{1/2}$. Hence, assumption (A.6) holds with $d_n = C\ell$, $n \geq 1$. $\square$

6.4. *Proofs of the results from Section 3.*

PROOF OF THEOREM 3.1. Let $\Sigma_{1n}$ be the $d \times d$ symmetric matrix with $(p,q)$th element given by $E(n^{-1}\sum_{i=1}^n Y_{in}^\#(p,q))$, $(p,q) \in \Lambda_0$ where $Y_{in}^\#$ is as in (2.1). Also, let

$$\hat{\Sigma}_n = \bigg[ \hat{\Gamma}_n(0) + \sum_{k=1}^\ell w_{kn}\{\hat{\Gamma}_n(k) + \hat{\Gamma}_n(k)'\} \bigg].$$

Then, it can be shown after some algebra that for all $(p,q) \in \Lambda_0$,

$$\hat{\Sigma}_n(p,q) - \Sigma_{1n}(p,q)$$

$$= n^{-1}\sum_{i=1}^n Y_{in}(p,q) + b^{-3/2}\xi_n(p,q)$$

$$- 2\bigg[ \sum_{k=0}^\ell (1 + n^{-1}k)w_{kn}\bar{X}_{n,p}\bar{X}_{n,q} \bigg]$$

(6.46)

$$= b_n^{-1/2}\bigg[ \frac{1}{\sqrt{n\ell}}\sum_{i=1}^n Y_{in}(p,q) \bigg]$$

$$+ b_n^{-1}\bigg( -2\ell^{-1}\sum_{k=0}^\ell (1 + n^{-1}k)w_{kn} \bigg) Z_n^{e_p+e_q} + b^{-3/2}\xi_n(p,q)$$

$$\equiv b_n^{-1/2}\hat{A}_{1n}(p,q) + b_n^{-1}\hat{A}_{2n}(p,q) + b_n^{-3/2}\hat{A}_{3n}(p,q), \qquad \text{say},$$

where $\bar{X}_{n,p}$ is the $p$th component of $\bar{X}_n$ and where $\xi_n(p,q)$ is the $(p,q)$th element of $\xi_n \equiv n^{1/2}\,\text{SVEC}(\bar{X}_n[_3W_{1n} + {_3W_{b_0n}}]')$, where ${_3W_{1n}} = \ell^{-3/2}\sum_{i=1}^{\ell}(\sum_{k=i}^{\ell} w_{kn})X_i$ and ${_3W_{b_0n}} = \ell^{-3/2}\sum_{i=1}^{\ell}(\sum_{k=i}^{\ell} w_{kn})X_{n-i+1}$. Set $\hat{A}_{kn}(p,q) = \hat{A}_{kn}(q,p)$ for $p > q$. Then, $\hat{A}_{kn}$ are $d \times d$ symmetric matrices with $\|\hat{A}_{kn}\| = O_p(1)$ for $k = 1, 2, 3$. Next, by Taylor's expansion, for any $p = 1, \ldots, d$,

$$(6.47) \quad e_p' h(\bar{X}_n) = D^{e_p}H(\bar{X}_n) + \sum_{|\alpha|=1}^{\nu(s)} n^{-|\alpha|/2} D^{\alpha+e_p}H(\mu) \cdot \frac{Z_n^{\alpha}}{\alpha!} + R_{1n},$$

where, on the set $\{\|\bar{X}_n - \mu\| \leq \delta_0\}$, the remainder term $R_{1n}$ satisfies $|R_{1n}| \leq \sup\{|D^{\alpha}H(x)| : \|x-\mu\| \leq \delta_0, |\alpha| = \nu(s)+2\}\|Z_n\|^{\nu(s)+1} n^{-[\nu(s)+1]/2}$, for some $\delta_0 > 0$ [determined by (A.1)(ii)]. In particular, $R_{1n} = \tilde{O}_p(\delta_{n,C})$.

Similarly, using assumptions (A.1)–(A.3) and moderate deviation inequalities for $\bar{X}_n$ [cf. Theorem 2.11, [GH], Theorem 2.4, Lahiri (1993)] and for sums of block variables $W_{jn}$'s (cf. Theorem 2.4, [L]), and Taylor's expansion, we get

$$(6.48) \quad \left. \begin{array}{l} n^{1/2}(\hat{\theta}_n - \theta) = \displaystyle\sum_{j=1}^{\nu(s)} n^{-(j-1)/2} \sum_{|\alpha|=j} \frac{D^{\alpha}H(\mu)}{\alpha!} Z_n^{\alpha} + \tilde{O}_p(\delta_{n,C}), \\ \hat{\tau}_n^{-1} = \tau_{1n}^{-1} + \displaystyle\sum_{j=1}^{s-2} \frac{g^{(j)}(\tau_{1n}^2)}{j!}(\hat{\tau}_n^2 - \tau_{1n}^2)^j + \tilde{O}_p(\delta_{n,C}), \end{array} \right\}$$

where $g(x) = x^{-1/2}$, $x > 0$, and $g^{(j)}(x) = \frac{(-1)^j 2j!}{2^{2j} j!} x^{-\frac{2j+1}{2}}$ is its $j$th derivative at $x \in (0, \infty)$, $j \geq 1$. Combining (6.46), (6.47) and (6.48), after some lengthy algebra, we obtain

$$T_n = T_{1n} + \tilde{O}_p(\delta_n),$$

where $T_{1n}$ is as in (3.1). Thus, the $(s-2)$th-order Edgeworth expansions of $T_n$ and $T_{1n}$ coincide, upto an error of order $O(\delta_n)$. It is clear that $T_{1n}$ is a smooth (infinitely differentiable) function of $\zeta_n \equiv (n^{1/2}\bar{X}_n', \frac{1}{\sqrt{n\ell}}\sum_{i=1}^{b_n} Y_{in}' + b_n^{-1}\tilde{\xi}_n')'$. Hence, using Theorem 5.1 and the transformation technique of Bhattacharya and Ghosh (1978), the $(s-2)$th-order EE for $T_{1n}$ can be derived from a $(s-2)$th-order EE for $\zeta_n$. This completes the proof of the theorem. $\square$

PROOF OF PROPOSITION 3.2. Write $Z_{jn} = \sum_{|\alpha|=j+1} D^{\alpha}H(\mu)Z_n^{\alpha}/(\alpha!)$, $j \geq 0$. Then, using (6.46) and (6.47), after some lengthy algebraic manipulations, one can show that

$$\hat{\tau}_{1n}^2 - \tau_{1n}^2 = b_n^{-1/2}Q_{1n} + n^{-1/2}Q_{2n} + n^{-1/2}b_n^{-1/2}Q_{3n}$$

(6.49)
$$+ b_n^{-1}Q_{4n} + n^{-1}Q_{5n} + b_n^{-3/2}Q_{6n} + \tilde{O}_p(\delta_{n,C}),$$

where

$$Q_{1n} = h(\mu)'\hat{A}_{1n}h(\mu), \qquad Q_{2n} = 2\sum_{|\alpha|=1} Z_n^\alpha [\alpha' H^{(2)}(\mu)]\Sigma_{1n}h(\mu),$$

$$Q_{3n} = h(\mu)'\hat{A}_{2n}h(\mu), \qquad Q_{4n} = 2\sum_{|\alpha|=1} Z_n^\alpha \alpha' H^{(2)}(\mu)\hat{A}_{1n}h(\mu),$$

$$Q_{5n} = \sum_{|\alpha|=1}\sum_{|\beta|=1} Z_n^{\alpha+\beta}[2H^{(3)}_{\alpha,\beta}(\cdot;\mu)'\Sigma_{1n}h(\mu) + \alpha' H^{(2)}(\mu)\Sigma_{1n}H^{(2)}(\mu)\beta],$$

$$Q_{6n} = h(\mu)'\hat{A}_{3n}h(\mu).$$

Here $H^{(2)}(\mu) = ((D^{e_i+e_j}H(\mu)))$ is the $d \times d$ matrix of second-order partial derivatives of $H$ at $\mu$, and $H^{(3)}_{\alpha,\beta}(\cdot;\mu)$ is the $d \times 1$ vector with $p$th component $D^{e_i+e_j+e_p}H(\mu)$, $p = 1, \ldots, d$. Then stochastic approximation $T_{1n}$ is now given by

(6.50)
$$T_{1n} = \sum_{j=1}^{7} a_{jn}\hat{G}_{jn},$$

where $a_{1n} = 1$, $a_{2n} = b_n^{-1/2}$, $a_{3n} = n^{-1/2}$, $a_{4n} = b_n^{-1}$, $a_{5n} = b_n^{-1/2}n^{-1/2}$, $a_{6n} = n^{-1}$, $a_{7n} = b_n^{-3/2}$, and where with $\gamma_{jn} = g^{(j)}(\tau_{1n}^2)/j!$, $j \geq 1$,

$$\hat{G}_{1n} = \tau_{1n}^{-1}Z_{0n},$$
$$\hat{G}_{2n} = \gamma_{1n}Z_{0n}Q_{1n},$$
$$\hat{G}_{3n} = \tau_{1n}^{-1}Z_{1n} + \gamma_{1n}Z_{0n}Q_{2n},$$
$$\hat{G}_{4n} = \gamma_{1n}Z_{0n}Q_{3n} + \gamma_{2n}Z_{0n}Q_{1n}^2,$$
$$\hat{G}_{5n} = \gamma_{1n}Z_{0n}Q_{4n} + \gamma_{1n}Z_{1n}Q_{1n} + 2\gamma_{2n}Z_{0n}Q_{1n}Q_{2n},$$
$$\hat{G}_{6n} = \tau_{1n}^{-1}Z_{2n} + \gamma_{1n}Z_{1n}Q_{2n} + \gamma_{2n}Z_{0n}Q_{2n}^2 + \gamma_{1n}Z_{0n}Q_{5n},$$
$$\hat{G}_{7n} = \gamma_{1n}Z_{0n}Q_{6n} + \gamma_{3n}Z_{0n}Q_{1n}^3 + 2\gamma_{2n}Z_{0n}Q_{1n}Q_{3n}.$$

Note that by arguments in the proof of Lemma 3.28 of [GH], for unit vectors $x \in \mathbb{R}^{d_1}$ and $y \in \mathbb{R}^{d_1}$,

$$|E(x'\hat{A}_{1n})(y'Z_n)|$$
$$= b_n^{-1}\left|\sum_{i=1}^{b_{0n}}\sum_{j=1}^{b_{0n}} \operatorname{Cov}(x'[_2W_{in}], y'[_1W_{jn}])\right|$$

$$(6.51) \quad \leq b_n^{-1} \left| \sum_{i=1}^{b_{0n}} \sum_{|j-i| \leq C} \text{Cov}(x'[_2 W_{in}], y'[_1 W_{jn}]) \right| + O(\exp(-C(\kappa)\ell))$$

$$\leq C w_n \max_{1 \leq p,q,r \leq d} \ell^{-3/2} \sum_{i=1}^{C\ell} \sum_{j=1}^{C\ell} \sum_{k=1}^{C\ell} |E(e'_p X_i, e'_q X_j, e'_r X_k)|$$

$$+ O(\exp(-C(\kappa)\ell))$$

$$= O(\ell^{-1/2}),$$

provided $w_n \equiv \max\{|w_{kn}| : 1 \leq k \leq \ell\} = O(1)$ as $n \to \infty$. Thus, $\ell^{1/2} E\hat{G}_{2n} = O(1)$. By similar arguments, $|E\hat{G}_{4n}| + |E\hat{G}_{7n}| = O(b_n^{-3/2}\ell^{-1/2})$. Next note that by (6.13) above and by Lemma 4.1 of [L], $E\hat{G}_{5n} = O(b_n^{-1/2})$. Also, by Lemma 3.28 of [GH], $E\hat{G}_{6n} = O(n^{-1/2})$. Hence,

$$\mathcal{X}_{1,n} = \sum_{j=1}^{7} a_{jn} E(\hat{G}_{jn})$$

$$(6.52) \quad = n^{-1/2}[\ell^{1/2} E\hat{G}_{2n} + E\hat{G}_{3n}] + O(b_n^{-1} n^{-1/2})$$

$$\equiv n^{-1/2} \beta_{1,1,n} + O(b_n^{-1} n^{-1/2}), \quad \text{say.}$$

Next, define the $p$th cumulant of random variables $V_1, \ldots, V_p$ by

$$\mathcal{K}_p(V_1, \ldots, V_p) = (-\iota)^r \frac{\partial}{\partial \varepsilon_1} \cdots \frac{\partial}{\partial \varepsilon_p} \bigg|_{\varepsilon_1 = \cdots = \varepsilon_p = 0}$$

$$(6.53) \quad \times \log E \exp(\iota[\varepsilon_1 V_1 + \cdots + \varepsilon_p V_p])$$

(assuming that the relevant partial derivatives exist). Also, for $p_1 + \cdots + p_r = p$, $p_1, \ldots, p_r \in \mathbb{N}$ and random variables $W_1, \ldots, W_r$, write $\mathcal{K}_p(W_1^{\diamond p_1}, \ldots, W_r^{\diamond p_r}) = \mathcal{K}_p(W_1, \ldots, W_1, \ldots, W_r, \ldots, W_r)$ where on the right-hand side, $W_1$ appears $p_1$-times, $W_2$ appears $p_2$-times, etc. Then, using the well-known formula for expressing cumulants of polynomials of random variables [cf. Section 2.3, Brillinger (1981)] in terms of cumulants of the random variables themselves, Lemma 4.1 of [L], and arguments similar to (6.51), we have:

$$\mathcal{X}_{2,n} = \mathcal{K}_2(\hat{G}_{1n}^{\diamond 2})$$

$$+ b_n^{-1}[2b_n^{1/2}\mathcal{K}_2(\hat{G}_{1n}, \hat{G}_{2n}) + 2\mathcal{K}_2(\hat{G}_{1n}, \hat{G}_{4n}) + \mathcal{K}_2(\hat{G}_{2n}^{\diamond 2})]$$

$$+ n^{-1}[2n^{1/2}\mathcal{K}_2(\hat{G}_{1n}, \hat{G}_{3n})$$

$$(6.54) \quad + 2\ell^{1/2}\{\mathcal{K}_2(\hat{G}_{2n}, \hat{G}_{3n}) + \mathcal{K}_2(\hat{G}_{1n}, \hat{G}_{5n})\}$$

$$+ 2\mathcal{K}_2(\hat{G}_{1n}, \hat{G}_{6n}) + \mathcal{K}_2(\hat{G}_{3n}^{\diamond 2})]$$

$$+ O(b_n^{-2} + b_n^{-1} n^{-1/2})$$
$$\equiv e_n^2 + b_n^{-1} \beta_{2,1,n} + n^{-1} \beta_{2,2,n} + O(b_n^{-2} + b_n^{-1} n^{-1/2}), \quad \text{say};$$

$$\mathcal{X}_{3,n} = n^{-1/2}(n^{1/2}\mathcal{K}_3(\hat{G}_{1n}^{\diamond 3}) + 3\ell^{1/2}\mathcal{K}_3(\hat{G}_{1n}^{\diamond 2}, \hat{G}_{2n}) + 3\mathcal{K}_3(\hat{G}_{1n}^{\diamond 2}, \hat{G}_{3n}))$$
(6.55)
$$+ b_n^{-3/2}\mathcal{K}_3(\hat{G}_{1n}^{\diamond 2}, \hat{G}_{7n}) + O(b_n^{-1} n^{-1/2})$$
$$\equiv n^{-1/2} \beta_{3,1,n} + b_n^{-3/2} \beta_{3,2,n} + O(b_n^{-1} n^{-1/2}), \quad \text{say};$$

$$\mathcal{X}_{4,n} = b_n^{-1}[4b_n^{1/2}\mathcal{K}_4(\hat{G}_{1n}^{\diamond 3}, \hat{G}_{2n}) + 6\mathcal{K}_4(\hat{G}_{1n}^{\diamond 2}, \hat{G}_{2n}^{\diamond 2})]$$
$$+ n^{-1}[\mathcal{K}_4(\hat{G}_{1n}^{\diamond 4}) + 4n^{1/2}\mathcal{K}_4(\hat{G}_{1n}^{\diamond 3}, \hat{G}_{3n})$$
$$+ \ell^{1/2}\{12\mathcal{K}_4(\hat{G}_{1n}^{\diamond 2}, \hat{G}_{2n}, \hat{G}_{3n}) + 4\mathcal{K}_4(\hat{G}_{1n}^{\diamond 3}, \hat{G}_{5n})\}$$
(6.56)
$$+ 4\ell\mathcal{K}_4(\hat{G}_{1n}^{\diamond 3}, \hat{G}_{4n}) + 6\mathcal{K}_4(\hat{G}_{1n}^{\diamond 2}, \hat{G}_{3n}^{\diamond 2}) + 4\mathcal{K}_4(\hat{G}_{1n}^{\diamond 3}, \hat{G}_{6n})]$$
$$+ b_n^{-3/2} 4\mathcal{K}_4(\hat{G}_{1n}^{\diamond 3}, \hat{G}_{7n}) + O(b_n^{-2} + n^{-1} b_n^{-1/2})$$
$$\equiv b_n^{-1} \beta_{4,1,n} + n^{-1} \beta_{4,2,n} + b_n^{-3/2} \beta_{4,3,n} + O(b_n^{-1} n^{-1/2}), \quad \text{say};$$

where recall that $e_n^2 = \tau_n^2/\tau_{1n}^2$ and where the $\beta_{r,j,n}$-terms in (6.52)–(6.56) are $O(1)$ as $n \to \infty$. Further, by similar arguments, $\mathcal{X}_{5,n} = O(b_n^{-1} n^{-1/2})$.

Next, combining (6.52)–(6.56), and using the relation

$$\int_{\mathbb{R}} \exp(\iota t x)[\sigma^{-k} H_k(x/\sigma) \phi_\sigma(x)] \, dx = (\iota t)^k \exp(-t^2 \sigma^2/2), \quad t \in \mathbb{R},$$

one can express that the density of the preliminary EE $\psi_n^*$ for $T_n$ with error of order $o(n^{-1})$ as

$$\psi_n^*(x) = \phi_{e_n}(x)[1 + n^{-1/2} q_{1n}(x) + b_n^{-1} q_{2n}(x) + n^{-1} q_{3n}(x) + b_n^{-3/2} q_{4n}(x)],$$
$$x \in \mathbb{R},$$

where, with $H_{kn}(x) \equiv e_n^{-k} H_k(x/e_k), x \in \mathbb{R}, k \geq 1$,

$$q_{1n}(x) = \beta_{1,1,n} H_{1n}(x) + \beta_{3,1,n} H_{3n}(x)/6,$$
$$q_{2n}(x) = \beta_{2,1,n} H_{2n}(x)/2 + \beta_{4,1,n} H_{4n}(x)/24,$$
$$q_{3n}(x) = [\beta_{2,2,n} H_{2n}(x)/2 + \beta_{4,2,n} H_{4n}(x)/24]$$
(6.57)
$$+ [\beta_{1,1,n}^2 H_{2n}(x)/2 + \beta_{3,1,n}^2 H_{6n}(x)/72$$
$$+ \beta_{1,1,n} \beta_{3,1,n} H_{4n}(x)/6],$$
$$q_{4n}(x) = \beta_{3,2,n} H_{3n}(x)/6 + \beta_{4,3,n} H_{4n}(x)/24,$$

$x \in \mathbb{R}$.

Next, suppose that $a_n^{-1} \equiv (e_n^{-1} - 1) = O(n^{-1/3})$. Then using the relation

$$\int_{-\infty}^{x} H_{(k+1)n}(y)\phi_{e_n}(y)\,dy = -H_{kn}(x)\phi_{e_n}(x), \qquad x \in \mathbb{R},$$

and Taylor's expansions, one can expand the preliminary EE for $T_n$ to derive (3.8), where the respective polynomials $p_{in}(x)$, $x \in \mathbb{R}$, are given by

(6.58)
$$\begin{aligned}
p_{1n}(x) &= -[\beta_{1,1,n} + \beta_{3,1,n}e_n^{-3}H_2(x)], \\
p_{2n}(x) &= \beta_{3,1,n}e_n^{-3}xH_3(x), \\
p_{3n}(x) &= -\left[\frac{\beta_{2,1,n}}{2e_n^2}H_1(x) + \frac{\beta_{4,1,n}}{4!e_n^4}H_3(x)\right], \\
p_{4n}(x) &= \left[\frac{\beta_{2,1,n}}{2e_n^2}xH_2(x) + \frac{\beta_{4,1,n}}{4!e_n^4}xH_4(x)\right], \\
p_{5n}(x) &= -\left[\frac{\beta_{2,2,n} + \beta_{1,1,n}^2}{2}H_1(x) + \frac{\beta_{4,2,n}}{24}H_3(x)\right. \\
&\qquad \left. + \frac{\beta_{3,1,n}^2}{72}H_5(x) + \frac{\beta_{1,1,n}\beta_{3,1,n}}{6}H_3(x)\right], \\
p_{4n}(x) &= -[\beta_{3,2,n}H_2(x)/6 + \beta_{4,3,n}H_3(x)/24]. \qquad \square
\end{aligned}$$

6.5. *Proofs of results from Section 4.*

PROOF OF THEOREM 4.1. The block resampling estimator of $\text{Cov}(Z_n)$ is location invariant, and hence, w.l.g., we again set $\mu = 0$ (for notational simplicity). Let $\Sigma_{1n}^{[2]} = EU_{11}U_{11}'$. It is easy to check that

(6.59)
$$\begin{aligned}
\hat{\Sigma}_n^{[2]} - \Sigma_{1n}^{[2]} &= N^{-1}\sum_{i=1}^{N}(U_{1i}U_{1i}' - EU_{1i}U_{1i}') \\
&\quad + [w_n^0]^2 \bar{X}_n \bar{X}_n'\left[1 - \frac{2n}{N}\right] \\
&\quad + b_n^{-3/2}[w_n^0](\mathcal{E}_n^{[2]}Z_n' + Z_n\mathcal{E}_n^{[2]'}) \\
&\equiv b_n^{-1/2}\hat{A}_{1n}^{[2]}(p,q) + b_n^{-1}\hat{A}_{2n}^{[2]}(p,q) \\
&\quad + b_n^{-3/2}\hat{A}_{3n}^{[2]}(p,q), \qquad \text{say},
\end{aligned}$$

where $w_n^0 = \sum_{k=1}^{\ell} w_{kn}$, and

$$\mathcal{E}_n^{[2]} = [n/N]\ell^{-3/2}\sum_{k=1}^{\ell} w_{kn}\left(\sum_{i=1}^{k-1} X_i + \sum_{i=N+k}^{n} X_i\right)$$

is the *edge-effect* term. Note that the three terms in (6.59) are exactly of the same form as those in (6.46). Hence, an $(s-2)$th-order EE for $T_n^{[2]}$ can be derived retracing the proof of Theorem 3.1. We omit the details. $\square$

PROOF OF THEOREM 4.2. Since there is no mean correction, it is easy to verify that $\hat{\Sigma}_n^{[3]} - EU_{11}U_{11}'$ has the same expansion as in (6.59), with $\hat{A}_{kn} = 0$ for $k = 2, 3$. Hence, an EE for $T_n^{[3]}$ can be derived *directly* from the EE results of [L], using the transformation technique of Bhattacharya and Ghosh (1978). $\square$

## REFERENCES

ATHREYA, K. B. and LAHIRI, S. N. (2006). *Measure Theory and Probability Theory*. Springer, New York. MR2247694

BHATTACHARYA, R. N. and GHOSH, J. K. (1978). On the validity of the formal Edgeworth expansion. *Ann. Statist.* **6** 434–451. MR0471142

BHATTACHARYA, R. N. (1985). Some recent results on Cramer–Edgeworth expansions with applications. In *Multivariate Analysis* **6** 57–75. North-Holland, Amsterdam. MR0822288

BHATTACHARYA, R. N. and RANGA RAO, R. (1986). *Normal Approximation and Asymptotic Expansions*. Robert E. Krieger Publishing Co. Inc., FL. MR0855460

BRADLEY, R. C. (1983). Approximation theorems for strongly mixing random variables. *Michigan Math. J.* **30** 69–81. MR0694930

BRILLINGER, D. (1981). *Time Series. Data Analysis and Theory*. Holden-Day Inc., Oakland, CA. MR0595684

BUSTOS, O. H. (1982). General $M$-estimates for contaminated $p$th-order autoregressive processes: Consistency and asymptotic normality. Robustness in autoregressive processes. *Z. Wahrsch. Verw. Gebiete* **59** 491–504. MR0656512

GÖTZE, F. and HIPP, C. (1983). Asymptotic expansions for sums of weakly dependent random vectors. *Z. Wahrsch. Verw. Gebiete* **64** 211–239. MR0714144

GÖTZE, F. and HIPP, C. (1994). Asymptotic distribution of statistics in time series. *Ann. Statist.* **22** 2062–2088. MR1329183

GÖTZE, F. and KÜNSCH, H. R. (1996). Second-order correctness of the blockwise bootstrap for stationary observations. *Ann. Statist.* **24** 1914–1933. MR1421154

HALL, P. (1992). *The Bootstrap and Edgeworth Expansion*. Springer, New York. MR1145237

HALL, P., HOROWITZ, J. L. and JING, B.-Y. (1995). On blocking rules for the bootstrap with dependent data. *Biometrika* **82** 561–574. MR1366282

IBRAGIMOV, I. A. and LINNIK, Y. V. (1971). *Independent and Stationary Sequences of Random Variables*. Wolters-Noordhoff Publishing, Groningen. MR0322926

KÜNSCH, H. R. (1989). The jackknife and the bootstrap for general stationary observations. *Ann. Statist.* **17** 1217–1261. MR1015147

LAHIRI, S. N. (1993). Refinements in asymptotic expansions for sums of weakly dependent random vectors. *Ann. Probab.* **21** 791–799. MR1217565

LAHIRI, S. N. (1996a). On Edgeworth expansion and moving block bootstrap for studentized $M$-estimators in multiple linear regression models. *J. Multivariate Anal.* **56** 42–59. MR1380180

LAHIRI, S. N. (1996b). Asymptotic expansions for sums of random vectors under polynomial mixing rates. *Sankhyā Ser. A* **58** 206–224. MR1662519

LAHIRI, S. N. (2003). *Resampling Methods for Dependent Data*. Springer, New York. MR2001447

LAHIRI, S. N. (2007). Asymptotic expansions for sums of block-variables under weak dependence. *Ann. Statist.* **35** 1324–1350. MR2341707

LIU, R. and SINGH, K. (1992). Moving blocks jackknife and bootstrap capture weak dependence. In *Exploring the Limits of Bootstrap* (R. Lepage and L. Billard, eds.) 225–248. Wiley, New York. MR1197787

PAPARODITIS, E. and POLITIS, D. N. (2001). Tapered block bootstrap. *Biometrika* **88** 1105–1119. MR1872222

PETROV, V. V. (1975). *Sums of Independent Random Variables*. Springer, Berlin. MR0388499

PRIESTLEY, M. B. (1981). *Spectral Analysis and Time Series*. Academic Press, New York.

POLITIS, D. and ROMANO, J. P. (1994). Large sample confidence regions based on subsamples under minimal assumptions. *Ann. Statist.* **22** 2031–2050. MR1329181

POLITIS, D. N. and ROMANO, J. P. (1995). Bias-corrected nonparametric spectral estimation. *J. Time Ser. Anal.* **16** 67–103. MR1323618

TIKHOMIROV, A. N. (1980). On the convergence rate in the central limit theorem for weakly dependent random variables. *Theory Probab. Appl.* **25** 790–809. MR0595140

VELASCO, C. and ROBINSON, P. M. (2001). Edgeworth expansions for spectral density estimates and studentized sample mean. *Econometric Theory* **17** 497–539. MR1841819

DEPARTMENT OF STATISTICS
TEXAS A & M UNIVERSITY
COLLEGE STATION, TEXAS 77843-3143
USA
E-MAIL: snlahiri@stat.tamu.edu